\documentclass{scrartcl} 
\usepackage{amsmath,amssymb,amsthm}

\newtheorem{lemma}{Lemma}
\newtheorem{theorem}[lemma]{Theorem}
\newtheorem{remark}[lemma]{Remark}
\newtheorem{assumption}[lemma]{Assumption}
\usepackage{graphicx}
\usepackage{booktabs}
\usepackage{algorithm}
\usepackage{algpseudocode}
\usepackage{hyperref}
\usepackage{numprint}
\usepackage{MnSymbol}
\def\beqa{\begin{eqnarray*}}
\def\eeqa{\end{eqnarray*}}
\def\beqal{\begin{eqnarray}}
\def\eeqal{\end{eqnarray}}

\def\bu{\boldsymbol u}

\def\bole{\boldsymbol e}
\def\bw{\boldsymbol w}
\def\bx{\boldsymbol x}

\def\bV{\boldsymbol V}

\def\bv{\boldsymbol v}
\def\bbf{\boldsymbol f}
\def\bw{\boldsymbol w}
\def\bn{\boldsymbol n}

\def\bbf{\textbf f}
\def\bchi{\boldsymbol\chi}
\def\boleta{\boldsymbol\eta}
\def\bxi{\boldsymbol\xi}

\def\div{\mbox{div}\,}

\def\vector#1{\left(\begin{array}{c}#1\end{array}\right)}

\begin{document}

\title{Error analysis of a pressure correction method with explicit time stepping}

\author{Utku Kaya \and Thomas Richter\thanks{Institute for Analysis and Numerics, Otto von Guericke University of Magdeburg, Germany. \url{utku.kaya@ovgu.de}, \url{thomas.richter@ovgu.de}}}

\maketitle

\begin{abstract}
      The pressure-correction method is a well established approach for simulating unsteady, incompressible fluids. It is well-known that implicit discretization of the time derivative in the momentum equation e.g. using a backward differentiation formula with explicit handling of the nonlinear term results in a conditionally stable method. In certain scenarios, employing explicit time integration in the momentum equation can be advantageous, as it avoids the need to solve for a system matrix involving each differential operator. Additionally, we will demonstrate that the fully discrete method can be expressed in the form of simple matrix-vector multiplications allowing for efficient implementation on modern and highly parallel acceleration hardware.
      Despite being a common practice in various commercial codes, there is currently no available literature on error analysis for this scenario. In this work, we conduct a theoretical analysis of both implicit and two explicit variants of the pressure-correction method in a fully discrete setting. We demonstrate to which extend the presented implicit and explicit methods exhibit conditional stability. Furthermore, we establish a Courant–Friedrichs–Lewy (CFL) type condition for the explicit scheme and show that the explicit variant demonstrate the same asymptotic behavior as the implicit variant when the CFL condition is satisfied.

\end{abstract}

\section{Introduction}\label{sec1}
Pressure-correction methods stand out as a preferred choice when approximating the Navier-Stokes equations, with their notable advantage lying in how they manage velocity-pressure coupling within the time-stepping scheme. Importantly, they avoid the emergence of a system matrix with a saddle-point structure, replacing it with only a series of linear equation systems that can be solved successively.
Originating in the late 60s \cite{Chorin1968, Temam1968}, pressure-correction methods have undergone continuous refinement and thorough study within the literature. For a comprehensive overview, we direct readers to \cite{Guermond2006}.

The primary objective of this paper is to provide an error analysis for the pressure-correction scheme, specifically focusing on cases where a first-order explicit time integration is applied in the momentum equation. In the explicit method under consideration, the differential operators act solely on solutions from the previous time-step, and the solution for the current time-step is obtained by inverting the mass matrix. Although this approach is commonly employed, there is a noticeable absence of error analysis in the literature for the fully discretized case using simple finite element spaces. Our theoretical analysis adheres to the framework established in \cite{Guermond1998}, carefully tracking constants and restating the error results presented in the referenced work. Furthermore, we extend the analysis to cover the explicit variant within the same framework. We also show the analysis of an explicit method where the nonlinear term is approximated, such that it can be represented using matrix-vector multiplications and no reassembly of the system matrix is necessary after initialization.

The paper is organized as follows: In Section \ref{sec2}, we introduce the Navier-Stokes equations and the notation used in this paper. Section \ref{sec3} describes the temporal and spatial discretizations. The error analysis for both implicit and explicit variants is explained in Section \ref{sec4}. We provide two academic examples in Section \ref{sec5} to validate our theoretical analysis numerically.
\section{Problem description and notation} \label{sec2}
Let $\Omega \in \mathbb{R}^d, \; d=2,\,3,$ be a bounded polyhedral domain with boundary $\Gamma$. We consider the incompressible time-dependent Navier-Stokes equations
\begin{align}
\partial_t \bu - \nu \Delta \bu + \div (\bu \otimes \bu) + \nabla p & \;=\; \bbf \quad \text{ in } (0,T] \times \Omega,\label{eq:NS1}\\
 \div \bu & \;=\; 0  \quad  \text{ in } (0,T] \times \Omega, \label{eq:NS2}\\
\bu & \;=\; \mathbf{0} \quad  \text{ on } (0,T] \times \Gamma \label{eq:NS3}\\
\bu(0, \cdot) & \;=\; \bu_0 \quad  \text{ in } \Omega \label{eq:NS4}
\end{align}
where $\bu: (0,T] \times \Omega \to \mathbb{R}^d$ and  $p: (0,T] \times \Omega \to \mathbb{R}$ denote the velocity and pressure fields, respectively. The viscosity $\nu > 0 $ and external sources $\mathbf{f}$ are given.
\subsection{Notation}
We consider the usual Sobolev spaces $W^{m,p}(\Omega)$ with norm $\|\cdot\|_{m,p}$
and semi-norm  $|\cdot |_{m,p}$. In the case $p=2$, we set $H^n(\Omega)=W^{n,2}(\Omega)$. The $L^2$-inner product on $\Omega$ and the $L^2(\Omega)$-norm is denoted by $(\cdot,\cdot)_{\Omega}$ and $\|\cdot\|_0$, respectively. Finally, $\|\cdot\|_\infty$ is the essential supremum norm.
We further define 
\beqa
H^1_0(\Omega) \; &:=& \; \{ v \in H^1(\Omega) \;  : \; v|_{\Gamma}=0 \}, \\
L^2_0(\Omega) \;  &:=& \; \{ q \in L^2(\Omega) \;  : \; (q,1)=0 \}
\eeqa
and consider the function spaces \[\mathbf{V}:=[H^1_0(\Omega)]^d, \; Q:=L^2_0(\Omega)\text{ and } S:=H^1(\Omega).\]
The ansatz spaces $\mathbf{V}$ and $Q$ satisfy the following inf-sup condition
 \begin{align}
\exists \beta > 0 : \inf_{q\in Q\backslash\{0\} } \sup_{\bv \in \mathbf{V} \backslash\{\mathbf{0}\}} \frac{(q, \div \bv)}{\|q\|_0 |\bv |_1} \geq \beta .
\label{inf-sup}
\end{align}
We fix $M+1$ discrete time points $0=t_0<t_1<\ldots<t_M=T$ with a constant time step $k=t_n-t_{n-1}>0$, hence $t_n:=nk$ for $0\le n\le M=T/k$.
For any normed space equipped with the norm $\|\cdot\|_X$ and time interval $(t_j, t_n) \in (0,T]$ we use the semidiscrete function space $L^p(t_j, t_n;X) = \{ w^j, w^{j+1}, \dots, w^n  \in X\} $ with its norm 
\begin{align*}
 \| w \|_{L^p(t_j, t_n;X)} :=& \Big(\sum_{i=j}^n k \| w^i \|_X^p  \Big)^{\frac1{p}}, \; p\in [1,\infty),\\
  \| w \|_{L^{\infty}(t_j, t_n;X)} := & \max_{j \leq i \leq n} \| w^i \|_X^p. 
\end{align*}

In our theoretical analysis, we employ the notation $x \lesssim y$ to denote that there exists a positive constant $C>0$ satisying $x \leq C y.$ We will assume that this constant does not depend on the discretization parameters $h$ (in space) and $k$ (in time). When we use the Young's inequality
\[ xy \leq  \frac1{2\epsilon} x^2 + \frac{\epsilon}{2} y^2 \]
for arbitrary $\epsilon>0$, we write
$ xy  \lesssim c_Y x^2 + y^2 $
with $c_Y=\frac{1}{\epsilon^2}$ and the hidden constant $\frac{\epsilon}{2}$.
\section{Discretization}\label{sec3}
\subsection{Finite element spaces}
Let $\mathcal{T}_h(\Omega)$ be a shape-regular, admissible decomposition of $\Omega$ into 
quadrilaterals/hexahedra. We consider finite dimensional spaces $\bV_h \subset \bV, \; Q_h \subset Q$ and $S_h \subset S$ consisting of piecewise continous functions on $\mathcal{T}_h(\Omega).$
For the well-posedness of the discrete system, we require the discrete counterpart of the inf-sup condition \eqref{inf-sup}
to be valid, i.e.
 \begin{align}
\exists \beta_h\ge \beta > 0 : \inf_{q\in Q_h\backslash\{0\} } \sup_{\bv \in \mathbf{V}_h \backslash\{\mathbf{0}\}} \frac{(q, \div \bv)}{\|q\|_0 |\bv|_1} \geq \beta_h.
\label{disc_inf-sup}
\end{align}
Moreover, due to \eqref{disc_inf-sup} the space of weakly divergence-free functions is not empty, i.e.
\begin{align*}
    \bV_h^{div} := \{ \bchi \in \bV_h \; : \; (\div \bchi, \varphi) = 0 \quad \forall \varphi \in Q_h \} \neq \emptyset
    \end{align*}
\subsection{Temporal discretization of the pressure-correction method}
The method with backward Euler time discretization is as follows:
\begin{align} 
    \intertext{\textbf{Step 1:}  Find $\tilde \bu_h^n \in \bV_h$ such that}
    \frac1{k}(\tilde\bu_h^n - \bu_h^{n-1}, \bchi) + \nu(\nabla \tilde\bu_h^n, \nabla \bchi) + c( \tilde\bu_h^{n-1}, \tilde \bu_h^{n},\bchi) & \; = \; (\bbf^n,\bchi) +(p_{h}^{n-1},\div \bchi) \quad \forall \bchi \in \bV_h,
        \label{eq:pre_discrete_implemented}
    \intertext{where $c(\cdot,\cdot,\cdot)$ is defined as}
c( \bu, \bv,\bchi) &= -(\bu \otimes \bv, \nabla \bchi) -\frac1{2}((\div \bu) \bv, \bchi) \label{convective}
    \intertext{\textbf{Step 2:}  Find $p^n_h \in S_h$ such that}
    (\nabla (p^n_h - p^{n-1}_h), \nabla  \varphi)  & \; = \; -\frac1{k}(\div \tilde\bu_{h}^n, \varphi) \quad \forall \varphi \in S_{h}.\\
    \intertext{\textbf{Step 3:}  Find $\bu^n_h \in \bV_h$ such that}
    (\bu^n_h, \bchi)  & \; = \;(\tilde \bu^n_h, \bchi) + k(\nabla (p^n_h - p^{n-1}_h), \bchi) \quad \forall \bchi \in \bV_{h}.
    \label{eq:corr_discrete_implemented}
\end{align}
In \textbf{Step 1}, a predictor velocity is calculated, which is not necessarily divergence-free. Subsequently, in \textbf{Step 2}, a Poisson equation is solved to address the divergence error through pressure updates. Finally, in \textbf{Step 3}, the end-step velocity field is determined. An interesting property of this scheme is that while the predictor velocity $\tilde \bu$ adheres to the correct boundary conditions, the end-step velocity $\bu$ only satisfies the correct boundary conditions in the normal direction, i.e., 
$$\bu \cdot \bn = 0 \text{ but } \bu \times \bn \neq 0.$$ Moreover, the pressure $p$ has corner singularities, if the domain is not smooth. This scheme is known as the incremental pressure-correction scheme, as \textbf{Step 2} involves solving a Poisson equation to upgrade the pressure. In earlier literature, the original Chorin-Temam algorithm is referred to as the non-incremental form, and its time convergence is discussed in \cite{Rannacher2006}. For a comprehensive analysis of both time and space convergence of the incremental pressure-correction scheme, we follow \cite{Guermond1998}, where a framework for the fully discrete scheme is provided.

On the other hand, from an implementation standpoint, this approach results in a system matrix for \textbf{Step 1}, which needs updating at each time step due to the nonlinear term. The resulting linear equation system is then solved, for instance, using a geometric multigrid solver. To mitigate the computational costs associated with updating the system matrix at every time step and to find approximate solutions without solving a poorly conditioned linear system, we explore an explicit scheme, achieved by replacing \textbf{Step 1} with
$\textbf{Step 1}^{*}\textbf{:}$  Find $\tilde \bu_h^n \in \bV_h$ such that
\begin{align}
     \frac1{k}(\tilde\bu_h^n - \bu_h^{n-1}, \bchi) + \nu(\nabla \tilde\bu_h^{n-1}, \nabla \bchi) + c( \tilde\bu_h^{n-1}, \tilde \bu_h^{n-1},\bchi)  \; = \; (\bbf^{n},\bchi) +(p_{h}^{n-1},\div \bchi) \quad \forall \bchi \in \bV_h,
    \label{eq:ex_pre_discrete_implemented}
\end{align}
Both implicit and explicit variants are summarized in Algorithms \ref{alg:implicit} and \ref{alg:explicit}, respectively. The theoretical analysis of both methods are given in Sections \ref{sec4a} and \ref{sec4b}.

    \begin{algorithm}[t]
\begin{minipage}{0.46\textwidth}

        \centering
        \begin{algorithmic}[1]
            \State \text{Initialize pressure $p_h^0$, set $t=0,\; n=0$} 
            \For{$n \leq M$}
            \State $n \gets n+1$
            \State find $\tilde \bu^n$ from \textbf{Step 1} 
            \State find $p^n$ from \textbf{Step 2}
            \State find $\bu^n$ from \textbf{Step 3} 
            \EndFor
        \end{algorithmic}
    \end{minipage}
    \hfill
    \begin{minipage}{0.46\textwidth}
        \centering
        \begin{algorithmic}[1]
            \State \text{Initialize pressure $p_h^0$, set $t=0,\; n=0$} 
            \For{$n \leq M$}
            \State $n \gets n+1$
            \State find $\tilde \bu^n$ from $\textbf{Step 1}^*$ 
            \State find $p^n$ from \textbf{Step 2}
            \State find $\bu^n$ from \textbf{Step 3} 
            \EndFor
        \end{algorithmic}          
    \end{minipage}
    \caption{Implicit scheme}\label{alg:implicit}
            \caption{Explicit scheme}\label{alg:explicit}
            
    \end{algorithm}

The advantage of the fully explicit scheme is that the momentum equation can be approximated without the need to solve any linear systems. By an interpolation of the convective term into the finite element space, namely by approximating $\bu_h\otimes \bu_h$  with ${\cal I}_h(\bu_h\otimes \bu_h)$ we can further formulate the complete scheme in the form of matrix-vector products. That allows for highly efficient implementation using GPU accelerators based on pre-assembled matrices. The analysis in Section \ref{sec4c} provides the theoretical basis for this efficient algorithm.

\section{Error analysis}\label{sec4}
In this section we explore the error analysis of both the implicit and explicit versions of the pressure-correction method, as outlined in Algorithms \ref{alg:implicit} and \ref{alg:explicit}. While the analysis of the implicit method (Sec. \ref{sec4a}) primarily echoes the findings from \cite{Guermond1998}, the examination of the explicit method (Sec. \ref{sec4b}) is conducted to seamlessly integrate within the same framework.
Many of the error estimates are quite technical in nature. Therefore, we will first outline the path and move some of the evidence to the appendix. 

We will start in Section~\ref{est:prelim} with collecting some preliminary estimates and introducing notation required in the proofs. Section~\ref{sec4a} then deals with the implicit pressure correction scheme. Although this proof is given in literature we repeat it in detail as the further steps,  covering the explicit alternative, will closely follow it and vary in details only. The explicit part is split into Sections~\ref{sec4b} and \ref{sec4c}. While Section~\ref{sec4b} is the direct explicit counterpart of the implicit algorithm, we go one step further in Section~\ref{sec4c}. Here we will analyze an approximate reformulation of the convective term $\div(\bu\otimes \bu)$ that will allow us to base this term on a matrix-vector product to avoid costly assembly using numerical quadrature on a possibly unstructured mesh. Finally, Section~\ref{sec4d} will give an error estimate for the pressure in the explicit case.

\subsection{Notation and Preliminaries}\label{est:prelim}

We begin by providing estimates for the convective term, a crucial component in our error analysis.

\begin{lemma}[Estimating the nonlinearity]
  Following inequalities hold
    \begin{itemize}
\item     for each $(\bu ,  \bv, \bchi) \in   W^{2,2}(\Omega) \times H^1(\Omega) \times L^2(\Omega):$ 
    \begin{align}
        c( \bu, \bv,\bchi) \lesssim \Big( \|  \bu\|_{\infty}  + \| \bu\|_{1,3}\Big) \| \bv\|_{1,2} \|  \bchi\|_0 \label{nonlinest2},
    \end{align}
\item     for each $(\bu ,  \bv, \bchi) \in   W^{1,2}(\Omega) \times W^{2,2}(\Omega) \times L^2(\Omega):$ 
\begin{align}
    c( \bu, \bv,\bchi) \lesssim \|\bu\|_{1,2} \Big( \|  \bv\|_{\infty}  + \| \bv\|_{1,3}\Big)\| \bchi\|_0 \label{nonlinest3},
\end{align}
\end{itemize}
\end{lemma}
\begin{proof}
    The convective term is equivalent to
    \begin{align}
        c( \bu, \bv,\bchi) = ((\bu \cdot \nabla)\bv, \bchi) + \frac1{2}((\div \bu) \bv, \bchi)
    \label{convective-2}
    \end{align}
Using the Sobolev embedding $H^1(\Omega) \to L^6(\Omega)$ gives
    \begin{align*}
      c( \bu, \bv,\bchi) & \leq \|\bu\|_{\infty} \| \nabla \bv \|_0 \| \bchi \|_0 + \frac1{2} \| \div \bu\|_3 \| \bv\|_6 \|\bchi\|_0
      \lesssim \|\bu\|_{\infty}  \| \nabla \bv \|_0 \| \bchi \|_0 + \| \bu\|_{1,3} \| \bv\|_{1,2} \|\bchi\|_0
    \end{align*}
   hence \eqref{nonlinest2} follows. From \eqref{convective-2} we observe similarly
   \begin{align*}
     c( \bu, \bv,\bchi)  \leq &   \|\bu\|_{6}\|  \nabla \bv\|_3 \|\bchi\|_0  +  \|\nabla \bu\|_0\| \bv\|_{\infty}\| \bchi\|_0  
     \leq  \|\bu\|_{1,2} \Big( \|  \bv\|_{\infty}  + \| \bv\|_{1,3}\Big)\| \bchi\|_0.
    \end{align*}
\end{proof}
As next, we introduce the Gr\"onwall's Lemma in the form presented in \cite{Guermond1998}
\begin{lemma}[Grönwall's Lemma]
\label{lem:Gronwall}
Let $\delta,g_0, a_n,b_n,c_n, \gamma_n \geq 0$ for $n\in \{0,1,\dots\}$ such that
\begin{align*}
    a_n+ \delta \sum_{j=0}^n b_j \leq \delta \sum_{j=0}^n \gamma_j a_j + \delta \sum_{j=0}^n c_j + g_0.
\end{align*}
If $\gamma_j\delta \le 1$ for $ j \in \{0,1,\dots\}$, then we have
\begin{align*}
    a_n+ \delta \sum_{j=0}^n b_j \leq \exp \Big (\delta \sum_{j=0}^n \sigma_j\gamma_j\Big) \Big( \delta \sum_{j=0}^n c_j + g_0 \Big).
\end{align*}
for each  $n \leq 0$, with $\sigma_j = (1-\gamma_j\delta)^{-1}.$
\end{lemma}
Now, we introduce the elliptic projection of the velocity field onto the divergence free discrete solution space with its associated stability properties. Let $(\bw_h(t),r_h(t))\in \bV_h\times Q_h$ be the solution of
\begin{align}
\label{projection} 
\begin{split}
    \nu(\nabla \bw_h(t), \nabla \bchi) - (r_h, \div \bchi) &=  \nu(\nabla \bu(t),\nabla \bchi)  \qquad \forall \bchi  \in \bV_h,\\
    (\div \bw_h(t), \varphi) &= 0\qquad \forall \varphi \in Q_h,
\end{split}
\end{align}
For the solution of \eqref{projection} we have following stability and error estimate.
\begin{lemma}[Projection error]
\label{lem:stokesprojection}
Let $\bu^n\in \bV$ and $p^n\in Q$. Then, for the projection $\bw_h^n\in \bV_h$ and $r_h^n\in Q_h$ according to \eqref{projection} it holds
\[
\nu \|\nabla \bw^n_h\|_0 \leq   \nu \|\nabla\bu^n\|_0
\]
Furthermore, for $\bu^n\in \bV\cap H^{r+1}(\Omega)$ due to the inf-sup stability of the discrete solution spaces \eqref{disc_inf-sup} it holds
\begin{align}
\nu \| \nabla (\bu^n-\bw_h^n)  \|_0 + \beta_h\| r^n_h\|_0  &\leq  \inf_{\bv_h \in \bV_h} \nu\|\nabla(\bu^n - \bv_h) \|_0  \nonumber\\
& \lesssim   h^{r} \nu\|\bu^n\|_{r+1,2} \label{projection-error}
\end{align}
Moreover, the Aubin-Nitsche trick provides
\begin{align}
    \label{int_err_aubin_nitsche}
    \| \bu^n-\bw_h^n \|  &\lesssim  h^{r+1} \|\bu^n\|_{r+1,2} 
\end{align}
\end{lemma}
By standard Sobolev embeddings we have the following Lemma.
\begin{lemma}[Stability of the projection]\label{stab2}
If $\bu^n \in \bV \cap [W^{2,2}(\Omega)]^d$ and $d\le 3$ then it holds
\begin{align*}
    \| \bw_h^n \|_{ \infty}+ \| \bw_h^n \|_{1,3} \lesssim   \|\bu^n\|_{2,2} .
\end{align*}
\end{lemma}
We split the errors into interpolation errors $\beta$ and approximation errors $\bxi$ as follows
\begin{alignat*}{4}
\bole_u^n \; :=\;  \underbrace{u^n - \bw_h^n}_{=\boleta_u^n} + \underbrace{\bw_h^n - \bu_h^n}_{=\bxi_u^n},\qquad
\tilde \bole_u^n \; :=\;  \underbrace{u^n - \bw_h^n}_{=\boleta_u^n} + \underbrace{\bw_h^n - \tilde \bu_h^n}_{=\tilde\bxi_u^n},\qquad
e_p^n \; :=\; \underbrace{p^n - j_h p^n}_{=\eta_p^n} + \underbrace{j_h p^n - p_h^n}_{=\xi_p^n}.
\end{alignat*}  
Here $j_h: H^{r} \to Q_h$ is an interpolation operator satisfying the estimate
\begin{align}
\label{interpolation}
    \| p - j_h p \|_0 + h \|\nabla (p - j_h p) \|_0  +  \lesssim h^r \|p\|_{r,2} \qquad \forall p \in  H^{r}(\Omega).
\end{align}
We further set
\begin{align*}
    \delta_t \phi^n := \frac{1}{k}\big(\phi^n - \phi^{n-1}\big).
\end{align*}
We have for $\partial_t \phi^n \in {L^{\infty}(t_{n-1}, t_n; L^2(\Omega))}$
\begin{align}\label{time_derivative}
\|\delta_t \phi^n \|_0 \leq \| \partial_t \phi \|_{L^{\infty}(t_{n-1}, t_n; L^2(\Omega))}.
\end{align} 
moreover, for $\partial_{tt} \bu \in L^{\infty}(t_{n-1}, t_n;L^2(\Omega))$ we have
\begin{align}\label{time_approx}
  \| \partial_t \bu^n - \delta_t \bu^n \|_0 \leq \sqrt{k}\| \partial_{tt} \bu \|_{L^{2}(0,T;L^2(\Omega))}.
\end{align} 
and if $ \partial_{t} \bu^n \in H^{r+1}(\Omega)$, then we have
\begin{align}\label{time_approx_boleta}
  \| \delta_t \boleta_u^n \|_0 \lesssim h^{r+1} \| \delta_{t} \bu^n \|_{r+1,2} \lesssim \frac{h^{r+1}}{\sqrt{k}}\| \partial_{t} \bu \|_{L^{2}(0,T;W^{r+1,2}(\Omega))}.
\end{align} 


\subsection{Error analysis for the implicit pressure correction method}\label{sec4a}

To estimate the error of the implicit pressure correction scheme we analyze a reformulation of the algorithm in the following form.
Let $\mathbf{Y}_h:=\bV_h^{div} \oplus \nabla S_h$. The formulation \eqref{eq:pre_discrete_implemented} - \eqref{eq:corr_discrete_implemented} is equivalent to
\begin{align}
    \intertext{\textbf{Step 1:}  Find $\tilde \bu_h^n \in \bV_h$ such that}
    \frac1{k}(\tilde\bu_h^n - \bu_h^{n-1}, \bchi) + \nu(\nabla \tilde\bu_h^n, \nabla \bchi) + c( \tilde\bu_h^{n-1}, \tilde \bu_h^{n},\bchi) & \; = \; (\bbf^n,\bchi) +(p_{h}^{n-1},\div \bchi) \quad \forall \bchi \in \bV_h.
        \label{eq:pre_discrete}
    \intertext{\textbf{Step 2:}  Find $(\bu^n_h,p^n_h) \in \mathbf{Y}_h \times Q_h$ such that}
    \frac1{k}(\bu^n_h-\tilde\bu^n_h , \bchi) - (p^n_h - p^{n-1}_h, \div \bchi)  & \; = \;0 \quad \forall \bchi \in \bV_{h}, \label{eq:corr}\\
    (\div \bu^n_h, \varphi) & \; =\;  0 \quad \forall \varphi \in Q_{h} \nonumber
\end{align}
which we will use for our error analysis. For more details about this equivalence, we refer to \cite{Guermond1998}.

\begin{assumption}[Regularity assumption for the implicit error estimate]\label{ass:implicit}
We assume that the solution $\bu,p$ has the regularity
\[
\begin{aligned}
p&\in W^{1,\infty}(0,T;H^1(\Omega)),\\
\bu&\in H^1(0,T;H^{r+1}(\Omega))\cap H^2(0,T;L^2(\Omega)\cap W^{1,\infty}(0,T;H^1(\Omega))
\end{aligned}
\]
where $r$ is the spatial degree of the finite element approach. 
\end{assumption}

The first step of the error estimate is an inequality which tracks the error propagation of the approximation error $\bxi_u^n = \bw_h^n-\bu_h^n$ the (discrete) discrepancy between projection and solution. 
\begin{lemma}[Error propagation for the implicit pressure correction scheme]
\label{lem:1}
Let the solution $\bu,p$ satisfy Assumption~\ref{ass:implicit}.
 For each $n \in \{1, \dots, M\}$ we have
 \begin{multline*}
 \frac1{2}\big( \| \bxi_u^n \|_0^2 +  \| \tilde \bxi_u^n - \bxi_u^{n-1}\|_0^2 - \| \bxi_u^{n-1} \|_0^2 \big)  + (1-c_Y) k\nu \|\nabla \tilde \bxi_u^n\|_0^2 + \frac{k^2}{2}\| \nabla \xi_p^n\|_0^2 \\
  \lesssim \;  (C_{4,n}k + 4 c_Y)(\| \bxi_u^{n-1} \|_0^2 +  \| \tilde \bxi_u^n - \bxi_u^{n-1}\|_0^2) + c_Y(\| \bxi_u^{n-2} \|_0^2 +  \| \tilde \bxi_u^{n-1} - \bxi_u^{n-2}\|_0^2)\\
  + c_Y k\nu  \|\nabla \tilde \bxi_u^{n-1}\|_0^2+ \frac{k^2}{2}(1+k)\| \nabla \xi_p^{n-1}\|_0^2  +k\mathcal{F}^n
\end{multline*}
with 
\begin{multline*}
\mathcal{F}^n = \;  C_{1,n}k^3 + (C_{2,n}+C_{3,n}) h^{2(r+1)} 
+ \frac{k^2}{2}(k+1)  \|\partial_t p^n\|_{L^{\infty}(0, T; H^1(\Omega))}^2
+ h^{2r+2} \| \partial_{t}\bu^n \|_{L^{2}(0,T;W^{r+1,2}(\Omega))}^2\\
+ k^3  \| \partial_{tt}\bu^n \|_{L^{2}(0,T;L^2(\Omega))}^2
\end{multline*}
where $C_{1,n},\dots,C_{4,n}$ are independent of $h, k$ but may depend on the solution $(\bu, p)$ and the viscosity parameter $\nu$.
\end{lemma}
The lengthy proof is given in Appendix~\ref{app:implicit}. 

With Gronwäll's inequality we can estimate the resulting approximation error.
\begin{lemma}[Approximation error estimate for the implicit pressure correction scheme]
\label{theorem:1}
Let 
\[
C_n := (C_{4,n}+1)k +5 c_Y< 1,\quad C_{4,n} = 2\min \{(1 + c_p^2) \nu^{-1}, c_{inv}^2 k (1+h^{-2})\}   \|\bu^n\|_{2,2}^2
\]
for all $n \in \{1,\dots, N\}$. Then it holds 
    \begin{align*}
       \frac1{2} \big( \|  \bxi_u^m \|_0^2  + k^2\| \nabla \xi_p^m\|_0^2 \big) &+ \sum_{n=1}^m (1-2c_Y) k\nu \|\nabla \tilde \bxi_u^n\|_0^2 
      \; \lesssim \;   C_G k \sum_{n=1}^m\mathcal{F}^n
      \end{align*}
with 
\[
C_G \sim \exp\Big(\sum_{n=1}^{m}\frac{C_n}{1-C_n}\Big).
\]
\end{lemma}
\begin{proof}
Summing up the claim of Lemma \ref{lem:1} gives
\begin{multline*}
    \frac1{2} \big(\|  \bxi_u^m \|_0^2  + \| \tilde \bxi_u^m - \bxi_u^{m-1} \|_0^2 +k^2\| \nabla \xi_p^m\|_0^2) +  \sum_{n=1}^m (1-2c_Y)  k\nu \|\nabla \tilde \bxi_u^n\|_0^2 
\\
\lesssim 
\sum_{n=1}^m C_n \big( \|  \bxi_u^{n-1} \|_0^2 + \| \tilde \bxi_u^n - \bxi_u^{n-1} \|_0^2 + k^2 \| \nabla \xi_p^{n-1}\|_0^2\big) 
+ k \sum_{n=1}^m\mathcal{F}^n.
\end{multline*}
Finally, applying Gr\"onwall inequality, Lemma~\ref{lem:Gronwall}, gives the claim.
\end{proof}

The combination of the result with the interpolation error \eqref{int_err_aubin_nitsche} lets us estimate the error of the implicit pressure correction scheme.
\begin{theorem}[Error estimate for the implicit pressure correction scheme]
    \label{corollary:1}
Under the assumptions stated above and given sufficient regularity of $\bu,\; p$, it holds 
        \begin{align*}
    \| \bu - \bu_h \|_{L^{\infty}(0,T;L^2(\Omega) )}  +  \| \bu - \tilde \bu_h \|_{L^{\infty}(0,T;L^2(\Omega) )} & \; \lesssim \; h^{r+1} \| \bu \|^2_{L^\infty(0,T;W^{r+1,2})}   +  \Big( k C_G \sum_{n=1}^m\mathcal{F}^n \Big)^{\frac1{2}}\\
   \nu^{\frac1{2}} \| \nabla (\bu - \bu_h) \|_{L^{2}(0,T;L^2(\Omega) )}  & \; \lesssim \; h^{r} \| \bu \|_{L^2(0,T;W^{r+1,2})}   +  \Big( k C_G \sum_{n=1}^m\mathcal{F}^n \Big)^{\frac1{2}}
          \end{align*}
\end{theorem}

\subsection{Error analysis for the explicit pressure correction method}\label{sec4b}

We proceed with the explicit variant of the pressure correction scheme. The procedure is similar to the implicit part, but requires a different approach in some cases. To check step size conditions, it is necessary to determine constants precisely.

Like in the implicit part we start with a reformulation of Algorithm \ref{alg:explicit}, where the momentum step is treated explicitly:
\begin{align}
    \intertext{\textbf{Step 1:}  Find $\tilde \bu_h^n \in \bV_h$ such that}
    \frac1{k}(\tilde\bu_h^n - \bu_h^{n-1}, \bchi) + \nu(\nabla \tilde\bu_h^{n-1}, \nabla \bchi) + c( \tilde\bu_h^{n-1}, \tilde \bu_h^{n-1},\bchi) & \; = \; (\bbf^n,\bchi) +(p_{h}^{n-1},\div \bchi) \quad \forall \bchi \in \bV_h.
        \label{eq:explicit_pre_discrete}
    \intertext{\textbf{Step 2:}  Find $(\bu^n_h,p^n_h) \in \mathbf{Y}_h \times Q_h$ such that}
    \frac1{k}(\bu^n_h-\tilde\bu^n_h , \bchi) - (p^n_h - p^{n-1}_h, \div \bchi)  & \; = \;0 \quad \forall \bchi \in \bV_{h}, \label{eq:explicit_corr}\\
    (\div \bu^n_h, \varphi) & \; =\;  0 \quad \forall \varphi \in Q_{h}. \nonumber
\end{align}

The proof to the following lemma to estimate the discrete approximation error $\bxi_u^m$ differs only in parts from that of Lemma~\ref{lem:1} but mostly uses the same techniques. 

\begin{lemma}[Approximation error estimate for the explicit pressure correction scheme]
    \label{theorem_explicit}
    Let $\bu,p$ satisfy Assumption~\ref{ass:implicit} and further, let
   \[
   C_{*,n}  := C_n + k\big(\frac{\nu}{h^2}c^2_{inv} + \frac{\|\tilde \bu_h^{n-1}\|^2_{\infty}}{4\nu}\big) < 1,
   \]
   where $C_n$ is defined in Lemma~\ref{theorem:1}.
   Then it holds 
    \begin{align*}
   \frac1{2} \big(    \|  \bxi_u^m \|_0^2  + k^2\| \nabla \xi_p^m\|_0^2 \big) &+ \sum_{n=1}^m k\nu \|\nabla \tilde \bxi_u^n\|_0^2 
      \; \lesssim \;  k C_{G,*} \sum_{n=1}^m\mathcal{G}^n
      \end{align*}
with $C_{G,*} \sim \exp\Big(\sum_{n=1}^{m}\frac{C_{*,n} }{1-C_{*,n} }\Big)$ where  
\begin{equation}\label{GG}
\mathcal{G}^n:=\mathcal{F}^n + k^2\| p \|_{L^2(0,T;H^1(\Omega))}^2 + k(1+\nu k)\|\bu \|_{L^{\infty}(0,T;H^1(\Omega))}^2.
\end{equation}
    \end{lemma}
   \begin{proof}
We analyse the terms that contribute to the estimate in Theorem \ref{theorem:1}. Substracting \eqref{eq:explicit_pre_discrete} from \eqref{projection} gives
     \begin{align}\label{err-eq-explicit}
    \begin{split}
       (\tilde\bxi_u^n - \bxi_u^{n-1}, \bchi) + k\nu(\nabla \tilde\bxi_u^{n}, \nabla \bchi)  \; =& \;k(r_h^n - p_{h}^{n-1},\div \bchi)-  k(\bbf^n,\bchi) + (\bw_h^n - \bw_h^{n-1},\bchi) \\ 
       &+k\nu(\nabla \bu^n, \nabla \bchi) 
       -k\nu(\nabla (\tilde \bu_h^n-\tilde \bu_h^{n-1}), \nabla \bchi) + k c(\tilde \bu_h^{n-1} , \tilde  \bu_h^{n-1}, \bchi)\\
        =&- k(\partial_t \bu^n - \delta_t \bu^n, \bchi) - k( \delta_t \boleta_u^n, \bchi)
        -\underbrace{k\nu(\nabla (\tilde \bu_h^n-\tilde \bu_h^{n-1}), \nabla \bchi)}_{=:I} \\
       & - \underbrace{k \Big(c(\bu^{n} ,\bu^{n}, \bchi) -  c(\tilde \bu_h^{n-1} , \tilde  \bu_h^{n-1}, \bchi) \Big)}_{=:II} + k(j_h p^n - p_{h}^{n-1},\div \bchi)\\
       & +  k(r_h^n + p^n - j_h p^{n},\div \bchi)
    \end{split}
    \end{align}
Testing with $\bchi = \tilde\bxi_u^n$ gives  
\begin{align}
    I \leq& \sqrt{k\nu} h^{-1}c_{inv}\| \tilde \bu_h^n-\tilde \bu_h^{n-1} \|_0 \sqrt{k\nu} \|\nabla \tilde\bxi_u^n\|_0 \nonumber\\
    \leq& \sqrt{k\nu} h^{-1}c_{inv}(\| \tilde\bxi_u^n-\tilde \bxi_u^{n-1} \|_0 + k\| \delta_t \bw_h^n \|_0) \sqrt{k\nu} \|\nabla \tilde\bxi_u^n\|_0 \nonumber\\
    \lesssim & \frac{k\nu}{h^2}c^2_{inv}(\| \tilde\bxi_u^n-\tilde \bxi_u^{n-1} \|_0^2 + k^2\| \delta_t \bw_h^n \|_0^2) +   c_Y k\nu \|\nabla \tilde\bxi_u^n\|_0^2 \nonumber\\
    \lesssim&  \frac{k\nu}{h^2}c^2_{inv}\| \bxi_u^{n-1}-\tilde \bxi_u^{n-1} \|_0^2 + k\mathcal{G}^n  +   c_Y k\nu \|\nabla \tilde\bxi_u^n\|_0^2\label{est1}
\end{align}
and due to stability of \eqref{eq:corr} and $H^1$ stability of the projection $j_h$ we have
\begin{align*}
     \| \tilde\bxi_u^{n-1}- \bxi_u^{n-1} \|_0^2  & \leq   k^2  \|\nabla (p_h^{n-1} - p_h^{n-2}) \|_0^2\\
      & \leq  k^2 \Big( \|\nabla \xi_p^{n-1} \|_0^2 + \|\nabla \xi_p^{n-2} \|_0^2+ k^2\|\nabla \delta_t j_h p^{n-1} \|_0^2 \Big)\\
      & \leq  k^2 \Big( \|\nabla \xi_p^{n-1} \|_0^2 + \|\nabla \xi_p^{n-2} \|_0^2+ c k\| p \|_{L^2(0,T;H^1(\Omega))}^2 \Big).
\end{align*}
The nonlinear term is splitted as
\begin{align}
    II &= k c(\bu^{n} - \bu^{n-1} ,\bu^{n}, \bchi) + kc(\bu^{n-1} -  \bw_h^{n-1} ,\bu^{n}, \bchi) \nonumber\\
    &+kc( \bw_h^{n-1}, \bu^{n}- \bw_h^{n}, \bchi) +kc(  \bw_h^{n-1} - \tilde\bu_h^{n-1},\bw_h^{n}, \bchi) \nonumber\\
    & +kc(  \tilde \bu_h^{n-1} ,\bw_h^{n}-\tilde \bu_h^{n-1}, \bchi). \label{exp_nonlin_term}
\end{align}
Except last term, other terms are handled in Lemma \ref{lem:1}. For the last term we have by adding $\pm \bw_h^{n-1}$
\begin{align*}
    kc(  \tilde \bu_h^{n-1} ,\bw_h^{n}-\tilde \bu_h^{n-1}, \tilde\bxi_u^n)& = kc(  \tilde \bu_h^{n-1} ,k\delta_t\bw_h^n,\tilde\bxi_u^n) + kc(  \tilde \bu_h^{n-1} ,\tilde \bxi_u^{n-1} ,\tilde\bxi_u^n)
\end{align*}
By considering 
$$c( \bu, \bv,\bchi) = \frac1{2} \big(((\bu \cdot \nabla)\bv, \bchi) -( \bv \otimes \bu, \nabla \bchi) \big)$$
for the first estimate
\begin{align*}
kc(  \tilde \bu_h^{n-1} ,k\delta_t\bw_h^n,\tilde\bxi_u^n) \leq & \frac{k^2}{2} \|\tilde \bu_h^{n-1}\|_{\infty} \big( \|\nabla \delta_t\bw_h^n \|_0\| \tilde\bxi_u^n\|_0 +  \|\delta_t\bw_h^n \|_0 \|\nabla \tilde\bxi_u^n\|_0\big) \\
\lesssim & \frac{k^3}{4\nu} \|\tilde \bu_h^{n-1}\|^2_{\infty} \big( k\nu\|\nabla \delta_t\bw_h^n \|_0^2 +  \|\delta_t\bw_h^n \|^2_0 \big) + c_Y \| \tilde\bxi_u^n\|_0 + c_Y k\nu \|\nabla \tilde\bxi_u^n\|_0 \\
\lesssim & \frac{k^3}{4\nu} \|\tilde \bu_h^{n-1}\|^2_{\infty} \|\delta_t\bu^n \|_0^2 + c_Y \| \tilde\bxi_u^n\|_0 + c_Y k\nu \|\nabla \tilde\bxi_u^n\|_0 
\end{align*}
together with the $H^1$ stability of the projection \eqref{projection}. Moreover, we have
\begin{align*}
kc(  \tilde \bu_h^{n-1} ,\tilde \bxi_u^{n-1} ,\tilde\bxi_u^n) \leq & \frac{k}{2} \|\tilde \bu_h^{n-1}\|_{\infty}\big( \|\nabla \tilde \bxi_u^{n-1} \|_0\| \tilde\bxi_u^n\|_0 +  \|\tilde \bxi_u^{n-1} \|_0\|\nabla \tilde\bxi_u^n\|_0 \big) \\
\lesssim& \frac{k}{4 \nu}  \|\tilde \bu_h^{n-1}\|^2_{\infty}\big(\| \tilde \bxi_u^{n} \|^2_0+ \| \tilde \bxi_u^{n-1} \|^2_0 \big) + c_Y k\nu  \|\nabla \tilde \bxi_u^{n-1} \|^2_0 +   c_Y k\nu  \|\nabla \tilde \bxi_u^{n} \|^2_0.
\end{align*}
Using $\| \tilde\bxi_u^n \|_0^2 \leq \| \tilde\bxi_u^n- \bxi_u^{n-1} \|_0^2 + \| \bxi_u^{n-1} \|_0^2$
in the last estimate, the fact that $C_{*,n} < 1$ and applying Gr\"onwall inequality, Lemma~\ref{lem:Gronwall}, gives the claim.
\end{proof}

With this preparation, we can directly specify the a priori error estimate for the explicit algorithm.
\begin{theorem}[Error estimate for the explicit pressure correction scheme]
    \label{corollary:2}
Under the assumptions of Lemma~\ref{theorem_explicit}, in particular given that the time step restriction
\begin{equation}\label{expl:tsr}
  C_n + k\big(\frac{\nu}{h^2}c^2_{inv} + \frac{\|\tilde \bu_h^{n-1}\|^2_{\infty}}{4\nu}\big)< 1
\end{equation}
with $C_n$ given in Lemma~\ref{theorem:1}, it holds for the solutions of \eqref{eq:explicit_pre_discrete} - \eqref{eq:explicit_corr}
        \begin{align*}
    \| \bu - \bu_h \|_{L^{\infty}(0,T;L^2(\Omega) )}  +  \| \bu - \tilde \bu_h \|_{L^{\infty}(0,T;L^2(\Omega) )} & \; \lesssim \; h^{r+1} \| \bu \|^2_{L^\infty(0,T;W^{r+1,2})}   +  \Big( k C_{G,*} \sum_{n=1}^m\mathcal{G}^n \Big)^{\frac1{2}}\\
   \nu^{\frac1{2}} \| \nabla (\bu - \bu_h) \|_{L^{2}(0,T;L^2(\Omega) )}  & \; \lesssim \; h^{r} \| \bu \|_{L^2(0,T;W^{r+1,2})}   +  \Big( k C_{G,*} \sum_{n=1}^m\mathcal{G}^n \Big)^{\frac1{2}}
          \end{align*}
\end{theorem}

The error estimate shows the expected order of convergence. In contrast to the implicit pressure correction scheme the additional time step condition~\eqref{expl:tsr} must be satisfied.

\subsection{Error analysis for the explicit pressure correction method with approximated nonlinear term}\label{sec4c}

The explicit pressure correction method allows for very fast implementation and each time step mostly consists of several matrix vector products in the velocity space as well as the solution of the pressure Poisson problem that can be done efficiently by means of (geometric) multigrid method. One problematic term, however, is the non-linearity. On general grids and e.g. when using finite elements of higher order, it will be necessary to calculate this expression with numerical quadrature. Compared to matrix-vector products with sparse matrices, this is a considerably higher effort. This is particularly the case with unstructured or locally refined grids, where there is no clear structure and many indirect memory accesses are necessary for assembly. In the following, we propose an approximation that allows this term to be replaced with the help of matrix vector products. This approach takes ideas from mass lumping and is sometimes referred to as the \emph{fully practical finite element method}~\cite{Barrett2001}. 

Recall that functions $\bv_h \in \bV_h$ are of the form $\bv_h=\sum_i \bv_i \phi^h_i$ with $\bv_i \in \mathbb{R}^d$.
Now we show how the non-linearity can also be represented using matrix-vector multiplication applied to the nodal products $\bu_i\otimes\bu_j\in\mathbb{R}^{d\times d}$. In terms of finite element notation this discrete approximation means
\[
\mathcal{I}_h (\bv_h \otimes \bv_h)
=\sum_i (\bv_i\otimes \bv_i)\phi_i^h
\approx \Big(\sum_i \bv_i\phi_i^h\Big)\otimes\Big(\sum_j \bv_j\phi_j^h\Big)
\]
In the interpolation nodes it holds
\[
\mathcal{I}_h(\bv_h\otimes \bv_h)(\bx_k)=\bv_h(\bx_k) \otimes \bv_h(\bx_k).
\]
The benefit of this approximation is that it opens the possibility to implement the convective term as a matrix-vector product with a pre-assembled sparse matrix.
Note that, for a vector $\bv_j \in\mathbb{R}^3$, the product $\bv_j\otimes \bv_j\in\mathbb{R}^{3\times 3}$ is symmetric and hence it holds
\[
(\mathcal{I}_h(\bv_h\otimes \bv_h), \nabla \bchi^l_i)
=-\sum_j\sum_{k=1}^3 ((\bv_j\otimes  \bv_j)_{lk} \phi_j,\partial_k\phi_i) 
=-\sum_{k=1}^3\underbrace{\sum_j( \phi_j,\partial_k\phi_i)}_{=: M_k} (\bv_j\otimes \bv_j)_{lk} 
\]
where $\bchi_i^l = e_l\phi^h_i$ and $e_l \in\mathbb{R}^3$ is the unit vector in the coordinate direction $l$.
We therefore assemble and save the three sparse block matrices  $M_k = (m_{k,ij})_{i,j=1\dots Ndof}$ with $m_{k,ij} = (\phi^h_{j}, \partial_{k} \phi^h_{i})$, $k \in \{1,2,3\}$. Then, the nonlinear term is calculated as 
\[
n(\bv_h\otimes \bv_{h}, \bchi^l_i) = -\sum_{i=1}^{3} M_i (\bv^i\circ
\bv^l),\quad l=1,2,3,
\]
where we denote by $\bv^l=(\bv^l_1,\dots,\bv^l_{N_{dof}})$ the $l$-th
component of the velocity vector and by $(\bu\circ \bw)_i=\bu_i\bw_i$
for $i=1,\dots,N_{dof}$ the entry-wise product. 
Alltogether 6 such products  $\bv^1\circ \bv^1, \; \bv^1 \circ \bv^2,\; \bv^1 \circ \bv^3,\; \bv^2 \circ \bv^2,\; \bv^2 \circ \bv^3$
and $\bv^3 \circ \bv^3$ are calculated at each time step. 
Therefore, 
handling of the nonlinear term reduces to matrix-vector multiplications only and no reassembly of the residual is required.

This results in an algorithm that completely approximates the moment equation on the basis of sparse matrix-vector products and lumped mass matrices.
We give a brief analys for the error arising from the approximation of the nonlinear term considering the approximation 
\[
c^*( \tilde\bu_h^{n-1}, \tilde \bu_h^{n-1},\bchi) = -({\cal I}_h(\bu_h^{n-1}\otimes \bu_h^{n-1}), \nabla \bchi),
\]
in place of the nonlinear term in \eqref{eq:explicit_pre_discrete}. 

\begin{remark}
\label{rem1}
The claim of Lemma \ref{theorem_explicit} and therefore Corollary \ref{corollary:2} is still satisfied with 
\[\mathcal{\tilde G}^n = \mathcal{G}^n +  
4kc_{inv}^2c_P^2\|\nabla \tilde \bu_h^{n-1} \|_{\infty}^2 \| \bu \|_{L^{\infty}(0,T,H^1(\Omega))}^2  +  h^{2r} \nu\|\bu\|_{L^{\infty}(0,T,W^{r+1,2}(\Omega))}^2,
\]

where $\mathcal{G}^n$ is defined in~\eqref{GG}, if $C_{\ostar,n}<1$, where 
\[
C_{\ostar,n}  := C_{*,n} +  4k^2c^2_{inv}\|\nabla \tilde \bu_h^{n-1} \|_{\infty}^2.
\]
The only modification in the estimate of Lemma \ref{theorem_explicit} is at 
\begin{align*}
 k \Big(c(\bu^{n} ,\bu^{n}, \tilde \bxi_u^{n}) -  c^*(\tilde \bu_h^{n-1} , \tilde  \bu_h^{n-1}, \tilde \bxi_u^{n}) \Big)  = &  k \Big(c(\bu^{n} ,\bu^{n}, \tilde \bxi_u^{n}) -  c(\tilde \bu_h^{n-1} , \tilde  \bu_h^{n-1}, \tilde \bxi_u^{n}) \Big)  \\ &+  k \Big(c(\tilde \bu_h^{n-1} , \tilde  \bu_h^{n-1}, \tilde \bxi_u^{n})-  c^*(\tilde \bu_h^{n-1} , \tilde  \bu_h^{n-1}, \tilde \bxi_u^{n}) \Big).
\end{align*}
Note that the first term is exactly the term $II$ in the proof of Lemma \ref{theorem_explicit} and we only need to estimate the second term
\begin{align*}
    k \Big(c(\tilde \bu_h^{n-1} , \tilde  \bu_h^{n-1}, \tilde \bxi_u^{n})-  c^*(\tilde \bu_h^{n-1} , \tilde  \bu_h^{n-1},\tilde  \bxi_u^{n}) \Big)&  = k((\text{id}- {\cal I}_h)(\tilde \bu_h^{n-1}\otimes \tilde \bu_h^{n-1}), \nabla \tilde \bxi_u^{n})  -\frac{k}{2}((\div \bu_h^{n-1})\bu_h^{n-1},\tilde \bxi_u^{n}).
  \end{align*}
We observe using the interpolation property together with Young's inequality and the inverse inequality
\begin{align*}
k((\text{id}- {\cal I}_h)(\tilde \bu_h^{n-1}\otimes \tilde \bu_h^{n-1}),\nabla\tilde \bxi_u^{n}) & \leq  k\|(\text{id}- {\cal I}_h)(\tilde \bu_h^{n-1}\otimes \tilde \bu_h^{n-1})\|_0 \| \nabla \tilde \bxi_u^{n}\|_0
\lesssim k h\|\nabla(\tilde \bu_h^{n-1}\otimes\tilde  \bu_h^{n-1})\|_{0} \| \nabla \tilde \bxi_u^{n}\|_0\\
& \leq 2k h  \|\nabla \tilde \bu_{h}^{n-1} \|_{\infty}  \|  \tilde \bu_h^{n-1}\|_{0} \| \nabla \tilde \bxi_u^{n}\|_0\\
& \leq 2k c_{inv}  \|\nabla \tilde \bu_{h}^{n-1} \|_{\infty}  \|  \tilde \bu_h^{n-1}\|_{0} \| \tilde \bxi_u^{n}\|_0\\
& \leq  C_{\ostar,n}(\| \tilde \bxi_u^{n-1}\|^2_{0} + \| \bw_h^{n-1} \|^2_0  ) + c_Y \| \tilde \bxi_u^{n}\|^2_0
\end{align*}
as well as
\begin{align*}
    \| \bw_h^{n-1} \|^2_0  \lesssim c^2_P \|\nabla \bw_h^{n-1} \|^2_0 \lesssim c^2_P \| \nabla \bu^n \|_0^2.
\end{align*}
Moreover, we have 
\begin{align*}
\frac{k}{2}((\div \tilde \bu_h^{n-1})\tilde \bu_h^{n-1},\tilde \bxi_u^{n}) \leq & \frac{k}{2}\|\div \tilde \bu_h^{n-1} \|_0 \| \tilde \bu_h^{n-1}\|_{\infty} \|\tilde \bxi_u^{n} \|\\
\leq & \frac{k\nu}{4} \|\div \tilde \bu_h^{n-1} \|^2_0 + \frac{k}{4\nu}\| \tilde \bu_h^{n-1}\|^2_{\infty} \|\tilde \bxi_u^{n} \|^2_0  \\
\lesssim & \frac{k\nu}{4} (\|\nabla \tilde \bxi_u^{n-1} \|^2_0 + \|\nabla  \boleta_u^{n-1} \|^2_0 ) + \frac{k}{4\nu}\| \tilde \bu_h^{n-1}\|^2_{\infty} (\|\bxi_u^{n-1} \|^2_0 + \|\tilde\bxi_u^{n} -\bxi_u^{n-1}  \|^2_0 ),
\end{align*}
where the second term is already bounded, since the factor $\frac{k}{4\nu}\| \tilde \bu_h^{n-1}\|^2_{\infty}$ is already included in $C_{*,n}$.
\end{remark}
\subsection{Error analysis for the pressure field}\label{sec4d}
Now, we present the error analysis for the pressure field calculated from the explicit scheme \eqref{eq:explicit_pre_discrete}-\eqref{eq:explicit_corr}. We assume that the algorithm is initialized, such that the errors of the initial solutions vanish, i.e. for $\bxi_u=\bw_h-\bu_h$, $\tilde\bxi_u = \bw_h-\tilde \bu_h$ and $\xi_p=j_h p- p_h$  it holds
\begin{align}
\| \bxi_u^0 \|_0 = \|\tilde \bxi_u^0 \|_0 = \| \nabla \xi_p^0 \|_0 &= 0.
\label{assumption}
\end{align}

To formulate a pressure error estimate we must strengthen the assumptions.
\begin{assumption}[Regularity assumptions for the pressure a priori estimate]
  \label{ass:pressure}
In addition to Assumption~\ref{ass:implicit} let $\bu,p$ satisfy
\[
p\in L^2(0,T;H^r(\Omega))
\]
\end{assumption}

As the following proofs regarding the pressure error are even more technical we will refrain from tracking the constants in detail and will instead mostly use the notation $x\lesssim y$ with a generic constant that does not depend on $h>0$ and $k>0$. 

We start with an error estimate for the the very first time step that gives, compared to Lemma~\ref{theorem:1} and Lemma~\ref{theorem_explicit} a higher order approximation.
\begin{lemma}[Estimate for the initial approximation error]\label{approx:1:p}
Let Assumption~\ref{ass:pressure} hold. For the first time step of the explicit scheme \eqref{eq:explicit_pre_discrete}-\eqref{eq:explicit_corr} it holds 
\begin{align*}
\|\delta_t \tilde \bxi_u^1\|^2_0 + k \nu\|\delta_t \nabla \tilde \bxi_u^1\|^2_0 + k^2 \| \nabla \delta_t \xi_p^1\|_0^2  \lesssim   k^2 D_1 + h^{2r} D_2 
\end{align*}
with 
\[ 
    D_1 = \|\partial_{tt} \bu \|^2_{L^{\infty}(0,T;L^2(\Omega))} + \| \partial_t  p \|_{L^{\infty}(0,T;H^1(\Omega))}^2 +  (\| \bu^{0}  \|^2_{2,2}+\| \bu^{1}  \|^2_{2,2})\| \partial_t \bu\|^2_{L^{\infty}(0,T;H^1(\Omega))}
\]
and 
\[
D_2 = h^{2}\|\partial_t \bu^1\|_{L^{\infty}(0,T;H^{r+1}(\Omega))} + \|\bu^1\|^2_{r+1,2} + \|p\|_{r,2}^2 + \| \bu^0\|^2_{r+1,2}  \|\bu^1\|_{2,2}^2 + \| \bu^0\|^2_{3,2}  \|\bu^1\|_{r+1,2}^2
\]
\end{lemma}
The proof is a slight modification of the proof to Lemma~\ref{theorem_explicit} and we give it in Appendix~\ref{app:pressure}

With this preparation we proceed with an estimate of the (discrete) time derivative of the velocity's approximation error.
\begin{lemma}
Under the assumption above and assumptions of Lemma \ref{theorem_explicit}, for the explicit scheme \eqref{eq:explicit_pre_discrete}-\eqref{eq:explicit_corr} it holds
$$\|\delta_t \bxi_u \|_{L^{\infty}(0,T; L^2(\Omega))}  + \sqrt{\nu}\|\delta_t \tilde \bxi_u \|_{L^{2}(0,T; H^1(\Omega))} \lesssim  k  + h^r $$ 
\end{lemma}
\begin{proof}
The error equation for velocity increments is written as
\begin{align*}
    (\delta_t\tilde\bxi_u^n - \delta_t\bxi_u^{n-1}, \bchi) + k\nu(\nabla \delta_t\tilde\bxi_u^{n}, \nabla \bchi)  \; =
     & \; T_1^n(\bchi) + T^n_2(\bchi) +T^n_3(\bchi) + T^n_4(\bchi)+T^n_5(\bchi) + T^n_6(\bchi)
\end{align*}
we test with $\bchi = \delta_t\tilde\bxi_u^n$ and perform term by term estimates. The first three terms are estimated straightforward:
\begin{align*}
T_1^n(\delta_t\tilde\bxi_u^n) &= -k(\delta_t(\delta_t \bu^n - \partial_t \bu^n), \delta_t\tilde\bxi_u^n) \lesssim k^2 \| \delta_t(\delta_t \bu^n - \partial_t \bu^n)\|^2 + c_Y \| \delta_t\tilde\bxi_u^n \|^2 \\
& \lesssim k^2 \| \partial_{t} \bu\|_{L^{\infty}(0,T;L^2(\Omega))}^2 + c_Y \| \delta_t\tilde\bxi_u^n \|^2\\
T_2^n(\delta_t\tilde\bxi_u^n) &= -k(\delta_{tt} \boleta_u^n, \delta_t\tilde\bxi_u^n )
\lesssim k^2\| \delta_{tt} \boleta_u^n\|^2_0 + c_Y\|\delta_t\tilde\bxi_u^n \|_0^2\\
& \lesssim k^2 h^{2r+2}\| \partial_{tt} \bu\|^2_{L^{\infty}( 0,T;L^2(\Omega))} + c_Y\|\delta_t\tilde\bxi_u^n \|_0^2\\
T_3^n(\delta_t\tilde\bxi_u^n) &=k\nu(\nabla (\delta_t(\tilde\bxi_u^n-\tilde \bxi_u^{n-1})), \nabla \delta_t\tilde\bxi_u^n)
\leq k \nu h^{-1} c_{inv} \| \delta_t(\tilde \bxi_u^n-\tilde \bxi_u^{n-1})\|_0 \| \nabla \delta_t\tilde\bxi_u^n\|_0\\
& \lesssim  \| \delta_t(\tilde \bxi_u^n-\tilde \bxi_u^{n-1})\|_0^2  + c_Y k\nu\|\nabla \delta_t\tilde\bxi_u^n \|_0^2
\end{align*}
As estimating the fourth term is rather lengthy and technical, we separate this term as Lemma~\ref{lem:Appendix} given in Appendix~\ref{app:pressure}.
\begin{align*}
    T_4^n(\delta_t\tilde\bxi_u^n)& \lesssim  k^2 + h^{2r} +  c_Y\big( \| \delta_t\tilde\bxi_u^{n-1}\|_0^2 + k\nu\|\nabla \delta_t\tilde\bxi_u^{n-1}\|^2 + \| \delta_t\tilde\bxi_u^n\|_0^2 + k\nu\|\nabla \delta_t\tilde\bxi_u^n \|^2\big) \\
    & + \big(k \nu^{-1} \| \tilde\bxi_u^{n-1}\|_0^2 + k^2  \|\nabla  \tilde\bxi_u^{n-1}\|_0^2\big) \|\partial_t \bu\|_{L^{\infty(0,T;H^2(\Omega))}}^2
    \end{align*}
The fifth term gives the control over the presure-error increment and with the sixth term they both are obtained similarly to the proof of Lemma \ref{lem:1}.
 \begin{align*}
 -T_5^n(\delta_t\tilde\bxi_u^n) &= -k(\delta_t (j_h^p - p_h^{n-1}), \div \delta_t\tilde\bxi_u^n) \lesssim \frac{k^2}{2} \big( \| \nabla \delta_t \xi_p^{n}\|^2_0 + \| \nabla \delta_t \psi_h^{n}\|^2_0 + \| \nabla \delta_t (p_h - p_h^{n-1})\|^2_0 \big)\\
& \lesssim \frac{k^2}{2} \big( \| \nabla \delta_t \xi_p^{n}\|^2_0 + k^2(1 + \frac1{k})\| \nabla \delta_{tt} j_h p_h^{n}\|^2_0 + (1+k)\| \nabla \delta_{t} \xi_p^{n-1}\|^2_0 \big) +  \frac1{2}\| \delta_t \bxi_u^n \|^2_0 -  \frac1{2}\| \delta_t \tilde\bxi_u^n \|^2_0 \\
T_6^n(\delta_t\tilde\bxi_u^n) &= k(\delta_t (r_h^n +p^n - j_h^p) , \div \delta_t\tilde\bxi_u^n)  \lesssim  \frac{k}{\nu}\|\delta_t r_h^n\|_0^2+  c_{Y}k\nu \|\nabla \delta_t\bxi_u^n \|_0^2 + k^2\|\delta_t \boleta_p^n\|_0^2 + c_{Y}\|\delta_t\tilde \bxi_u^n \|_0^2\\
& \lesssim k h^{2r} (\| \delta_t\bu^n\|_{r+1,2}^2 + k\| \delta_t p^n\|_{r,2}^2) + c_{Y} k\nu\|\nabla \delta_t\bxi_u^n \|_0^2  + c_{Y}\|\delta_t\tilde \bxi_u^n \|_0^2
\end{align*} 
Applying the Gr\"onwall's Lemma, similar to the proof of Theorem \ref{theorem:1} gives the claim.
\end{proof}

\begin{theorem}[Pressure error for the explicit correction scheme]
  Under the assumptions stated above it holds 
    \begin{align*}
tle        \| p^n - p_h^n \|_{L^2(0,T;L^2(\Omega))} \lesssim k + h^r
ti    \end{align*}
\end{theorem}
\begin{proof}
Combining, the equations \eqref{eq:corr} and \eqref{error_equation} provide
\begin{align*}
    (\xi_p^{n}, \div \bchi) &= (j_h p^n - p^{n-1}, \div \bchi) + k^{-1}(\bxi_u ^n - \tilde \bxi_u^n, \div \bchi)\\
    & = k^{-1}(\bxi_u^n - \bxi_u^{n-1}, \bchi) + \nu(\nabla \tilde\bxi_u^{n}, \nabla \bchi) + (\delta_t \boleta_u^n, \bchi) + (\delta_t\bu^n - \partial_t \bu^n, \bchi) \nonumber \\
    &+\nu(\nabla (\tilde \bu_h^n-\tilde \bu_h^{n-1}),\nabla \bchi) + \Big(c(\bu^{n} ,\bu^{n}, \bchi) -  c(\tilde \bu_h^{n-1} , \tilde  \bu_h^{n-1}, \bchi) \Big)
    - (r_h^n + p^n - j_h p^{n},\div \bchi).
\end{align*}
Hence,the discrete inf-sup condition provides
\begin{align*}
    \beta_h\|\xi_p^{n}\|_0 & \leq \sup_{\bv \in \mathbf{V}_h \backslash\{\mathbf{0}\}} \frac{(\xi_p^{n}, \div \bv)}{ |\bv|_1} 
    \leq \| \delta_t \bxi_u^n\|_{-1}+\nu \|\nabla \tilde \bxi_u^{n}\| + \|\delta_t \boleta_u^n \|_{-1} + \| \delta_t\bu^n - \partial_t \bu^n\|_{-1}\\
   & + \sup_{\bv \in \mathbf{V}_h \backslash\{\mathbf{0}\}} \frac{\nu(\nabla (\tilde \bu_h^{n} - \tilde \bu_h^{n-1}),\nabla \bv)}{ |\bv|_1} + \sup_{\bv \in \mathbf{V}_h \backslash\{\mathbf{0}\}} \frac{c(\bu^{n} ,\bu^{n}, \bv) -  c(\tilde \bu_h^{n-1} , \tilde  \bu_h^{n-1}, \bv)}{ |\bv|_1} 
    + \| r_h^n \|_0 + \| \boleta_p^n\|_0. 
\end{align*}
We observe
\begin{align*}
\sup_{\bv \in \mathbf{V}_h \backslash\{\mathbf{0}\}} \frac{\nu(\nabla (\tilde \bu_h^{n} - \tilde \bu_h^{n-1}),\nabla \bv)}{ |\bv|_1} \leq \nu\|\nabla \bxi_u^n \|  +  \nu\|\nabla \bxi_u^{n-1} \| +  k\nu\|\nabla \delta_t\bw_h^{n} \| 
\end{align*}
From estimates in the proofs of Lemma \ref{lem:1} and Theorem \ref{theorem_explicit} one obtains
\begin{align*}
    \sup_{\bv \in \mathbf{V}_h \backslash\{\mathbf{0}\}} \frac{\Big(c(\bu^{n} ,\bu^{n}, \bv) -  c(\tilde \bu_h^{n-1} , \tilde  \bu_h^{n-1}, \bv) \Big)}{ |\bv|_1}  & \leq c_P(k\sqrt{C_{1,n}} +h^{r}\sqrt{C_{2,n}} + h^{r}\sqrt{ C_{3,n}}) \\
 & + \|\bxi_u^{n-1} \|_{1,2}\|\bu^n \|_{2,2} + k \|\delta_t \nabla \bw_h^n \|  + \|\bxi_u^{n-1} \|_{1,2}.
\end{align*}
The claim is shown by applying the Gr\"onwall's Lemma while estimating the terms depending on $\bxi_u$ by using the results from \ref{theorem_explicit}.
\end{proof}

\section{Numerical examples}\label{sec5}
In this section we compare the expected asymptotic error with the numerical ones. All experiments were performed with the finite element toolkit Gascoigne 3D~\cite{gascoigne}.

Let $\Omega:=(0,1)^3, \;I:=(0,T]$ and $\nu=10^{-3}.$ We choose $\bbf(\bx, t)$ such that Equation \ref{eq:NS1} is fulfilled for the exact solution
\beqa
  \bu(\bx,t) &=& \vector{ \sin(x+t)(\cos(z+t)-\sin(y+t))\\
                         -\cos(x+t)\cos(y+t) - \sin(y+t)\cos(z+t)\\
                         \sin(z+t)( \cos(y+t)-\cos(x+t))},\\
                          p(\bx,t) & =& \sin(x-y -z +t) +8\sin^3(\frac1{2})\sin(\frac1{2} - t).
\eeqa

\subsection{Temporal refinement}

\begin{figure}[t]   
  \centering
  \includegraphics[width=0.49\textwidth]{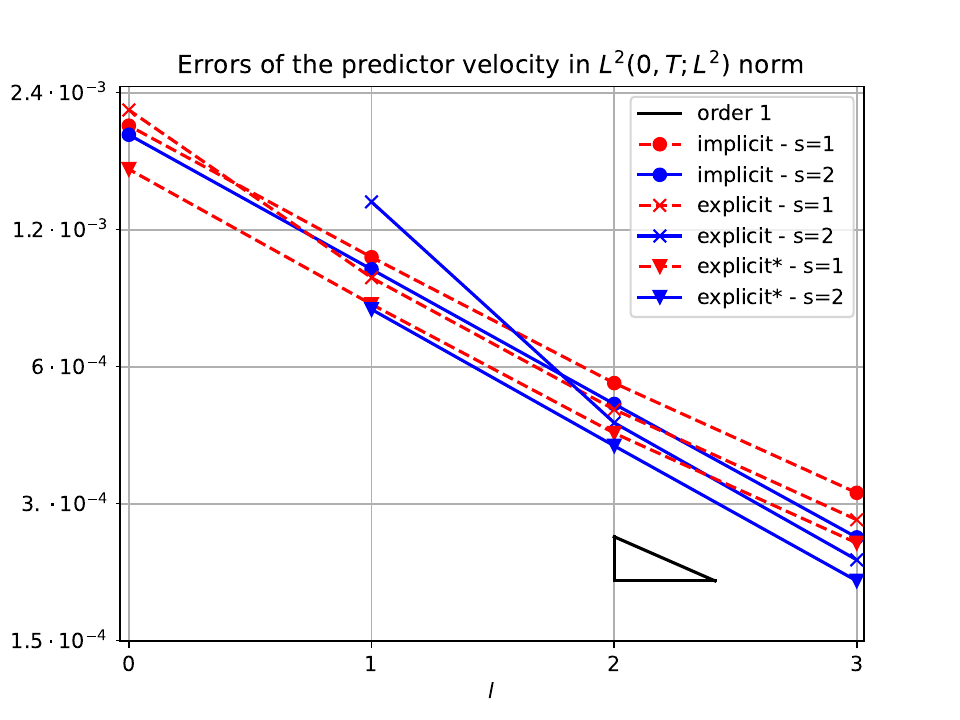}
  \includegraphics[width=0.49\textwidth]{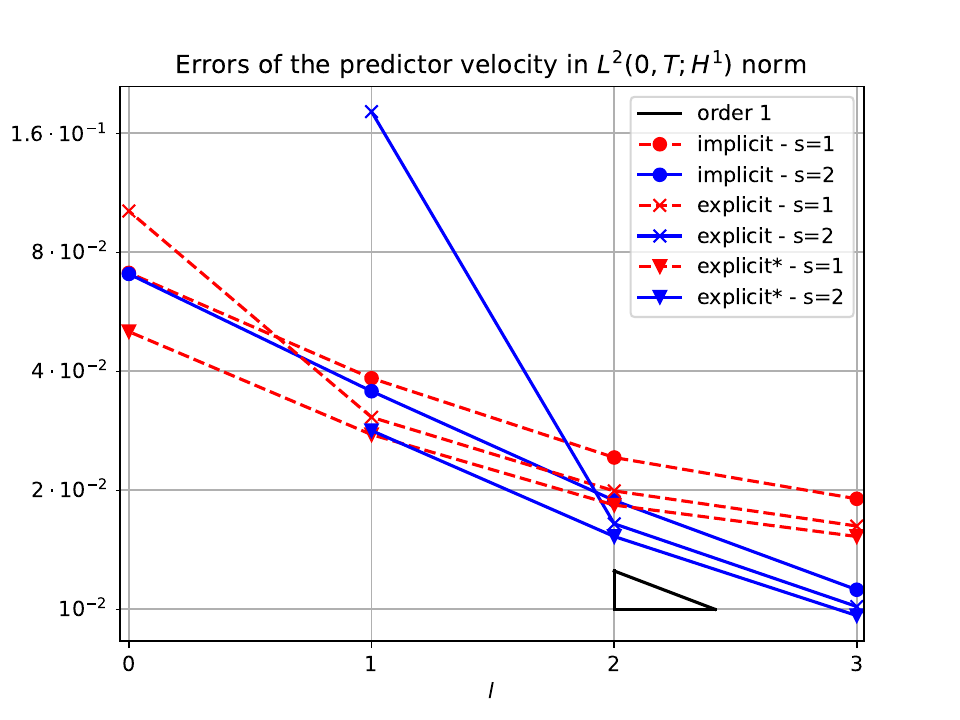}
  \caption{Errors of the predictor velocity in the $L^2(0,T;L^2)$-norm (left) and the $L^2(0,T;H^1)$-norm (right) under time-refinement (see Tables~\ref{tab:L2L2_k} and  \ref{tab:L2H1_k}).}
  \label{Fig1}
\end{figure}

We perform experiments with $T=1$ and choose the temporal and spatial discretization with the parameters
\[
\begin{aligned}
k &= 2^{-l-7},\; l \in \{ 0,\dots, 3\},\\
h &= 2^{-s-4},\; s \in \{1,2\},
\end{aligned}
\]
The errors of the predictor velocity are presented in Fig. \ref{Fig1}, both in  the spatial $L^2$-norm and the  $H^1$-norm. 
We observe from Figures \ref{Fig1} that, for all three presented methods; implicit, explicit and explicit* (the explicit method with approximated nonlinear term) for $k$ small enough the error in the $L^2(0,T;L^2(\Omega))$-norm converges to zero with each time refinement linearly as indicated with the theory. Similarly, in the $L^2(0,T;H^1(\Omega))$-norm (see Fig. \ref{Fig1} right) the error also decreases linearly in time. For coarse spatial meshes, $(s=1)$, full linear convergence is not reached for small time step sizes. 
Moreover, we observe from Figure~\ref{Fig1} that for both explicit and explicit* methods the choice $s=2,\; l <2$ does not exhibit linear convergence yet, as the CFL condition is not satisfied.

\begin{figure}[t]   
  \centering
  \includegraphics[width=0.49\textwidth]{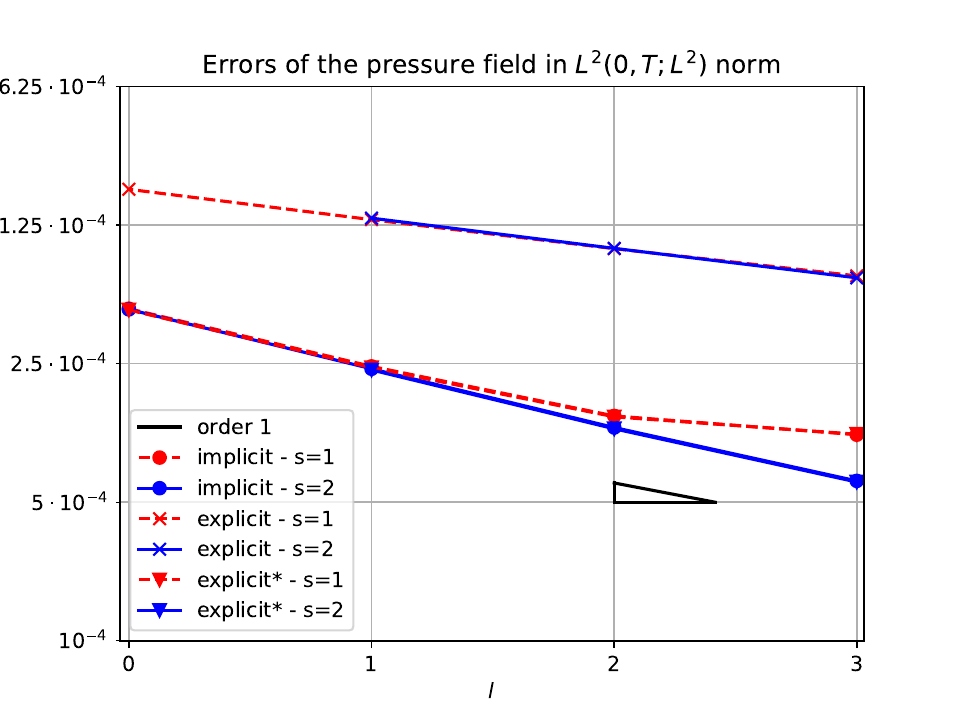}
  \caption{Errors of the pressure in the $L^2(0,T;L^2)$-norm under time-refinement (see Table \ref{tab:L2L2_pres_k})}
  \label{Fig3}
\end{figure}

\begin{table}[h]
\centering
\begin{tabular}{c|cc|cc|cc}
\toprule
$l$ & \multicolumn{2}{c|}{implicit}   & \multicolumn{2}{c|}{explicit} & \multicolumn{2}{c}{explicit*} \\
\midrule
 & $s=1$  & $s=2$ & $s=1$  & $s=2$ &$s=1$  & $s=2$ \\
\midrule
0 & $\numprint{2.03e-3}$ & $\numprint{1.94e-3}$ & $\numprint{2.20e-3}$ & $\numprint{1.30e1}$ & $\numprint{1.63e-3}$ & $\numprint{1.80e1}$ \\ 
1 & $\numprint{1.04e-3}$ & $\numprint{9.82e-4}$ & $\numprint{9.41e-4}$ & $\numprint{1.38e-3}$ & $\numprint{8.23e-4}$ & $\numprint{8.02e-4}$ \\ 
2 & $\numprint{5.52e-4}$ & $\numprint{4.97e-4}$ & $\numprint{4.83e-4}$ & $\numprint{4.52e-4}$ & $\numprint{4.29e-4}$ & $\numprint{4.02e-4}$ \\ 
3 & $\numprint{3.17e-4}$ & $\numprint{2.53e-4}$ & $\numprint{2.77e-4}$ & $\numprint{2.26e-4}$ & $\numprint{2.46e-4}$ & $\numprint{2.03e-4}$\\
\bottomrule
\end{tabular}
\caption{Errors of the predictor velocity in the $L^2(0,T;L^2)$-norm under temporal refinement}
\label{tab:L2L2_k}
\end{table}
\begin{table}[h]
\centering
\begin{tabular}{c|cc|cc|cc}
\toprule
$l$ & \multicolumn{2}{c|}{implicit}   & \multicolumn{2}{c|}{explicit} & \multicolumn{2}{c}{explicit*} \\
\midrule
 & $s=1$  & $s=2$ & $s=1$  & $s=2$ &$s=1$  & $s=2$ \\
\midrule
0 & $\numprint{7.09e-2}$ & $\numprint{7.05e-2}$ & $\numprint{1.02e-1}$ & $\numprint{7.75e0}$ & $\numprint{5.04e-2}$ & $\numprint{2.94e0}$ \\
1 & $\numprint{3.84e-2}$ & $\numprint{3.56e-2}$ & $\numprint{3.06e-2}$ & $\numprint{1.81e-1}$ & $\numprint{2.77e-2}$ & $\numprint{2.83e-2}$ \\
2 & $\numprint{2.42e-2}$ & $\numprint{1.88e-2}$ & $\numprint{1.99e-2}$ & $\numprint{1.65e-2}$ & $\numprint{1.83e-2}$ & $\numprint{1.53e-2}$ \\
3 & $\numprint{1.90e-2}$ & $\numprint{1.12e-2}$ & $\numprint{1.62e-2}$ & $\numprint{1.02e-2}$ & $\numprint{1.53e-2}$ & $\numprint{9.64e-3}$\\
\bottomrule
\end{tabular}
\caption{Errors of the predictor velocity in the $L^2(0,T;H^1)$-norm under temporal refinement}
\label{tab:L2H1_k}
\end{table}
\begin{table}[h]
\centering
\begin{tabular}{c|cc|cc|cc}
\toprule
$l$ & \multicolumn{2}{c|}{implicit}   & \multicolumn{2}{c|}{explicit} & \multicolumn{2}{c}{explicit*} \\
\midrule
 & $s=1$  & $s=2$ & $s=1$  & $s=2$ &$s=1$  & $s=2$ \\
\midrule
0 & $\numprint{4.73e-3}$ & $\numprint{4.69e-3}$ & $\numprint{1.89e-2}$ & $-$ & $\numprint{4.66e-3}$ & $-$ \\
1 & $\numprint{2.42e-3}$ & $\numprint{2.35e-3}$ & $\numprint{1.33e-2}$ & $\numprint{1.35e-2}$ & $\numprint{2.38e-3}$ & $\numprint{2.32e-3}$ \\
2 & $\numprint{1.36e-3}$ & $\numprint{1.19e-3}$ & $\numprint{9.49e-3}$ & $\numprint{9.53e-3}$ & $\numprint{1.35e-3}$ & $\numprint{1.17e-3}$ \\
3 & $\numprint{1.09e-3}$ & $\numprint{6.40e-4}$ & $\numprint{6.92e-3}$ & $\numprint{6.77e-3}$ & $\numprint{1.10e-3}$ & $\numprint{6.32e-4}$\\
\bottomrule
\end{tabular}
\caption{Errors of the pressure field in $L^2(0,T;L^2)$-norm under temporal refinement}
\label{tab:L2L2_pres_k}
\end{table}

\subsection{Spatial refinement}

\begin{figure}[t]   
  \centering
  \includegraphics[width=0.49\textwidth]{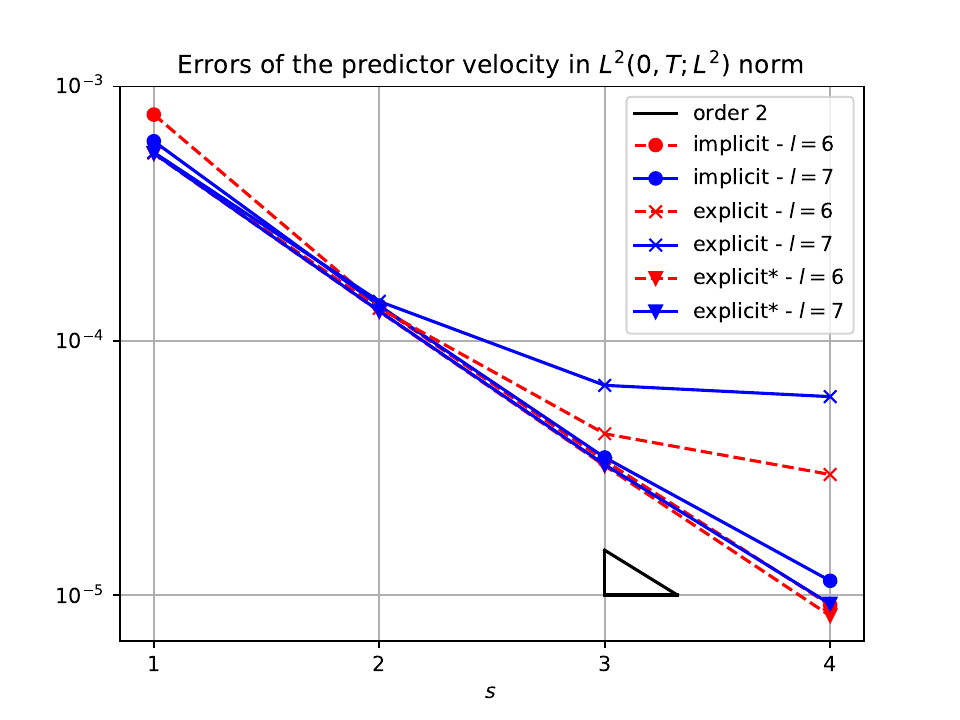}
  \includegraphics[width=0.49\textwidth]{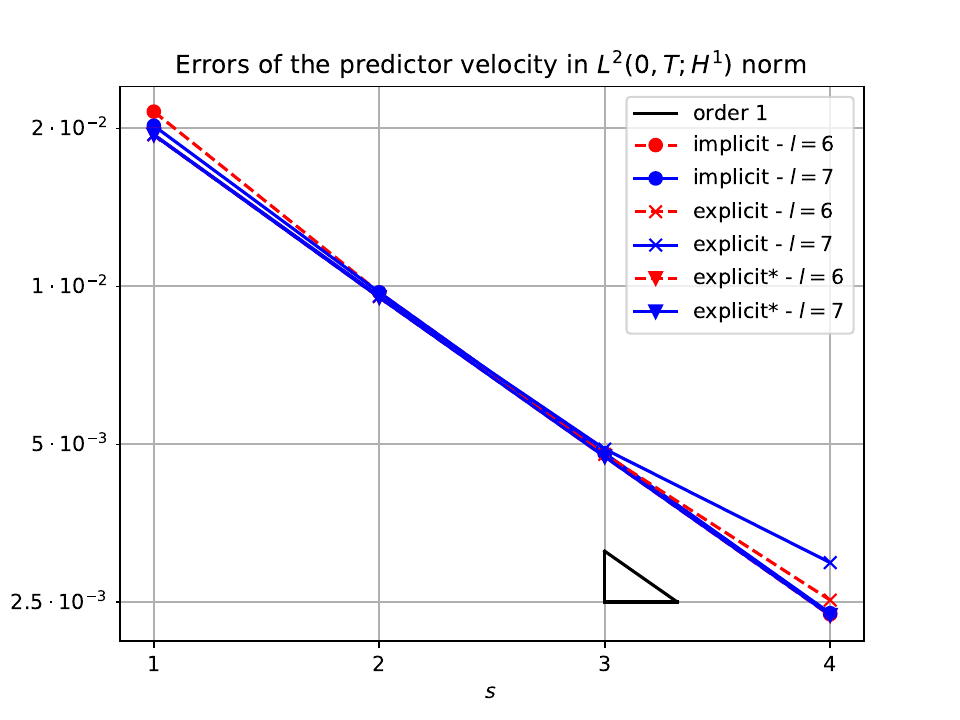}
  \caption{Errors of the predictor velocity in the $L^2(0,T;L^2)$-norm (left) and $L^2(0,T;H^1)$-norm (right) under spatial refinement (see Tables~\ref{tab:L2L2_h}) and~\ref{tab:L2H1_h}.}
  \label{Fig4}
\end{figure}

Next, we fix $T=0.1,$ and perform experiments with $k = 0.1 \cdot 2^{-l-7},\; l \in \{ 0,\dots, 1\}$ and $h = 2^{-s-4}$ for $s\in \{1,2,3,4\}$. For this test case, the errors of the predictor velocity are presented in Fig. \ref{Fig4}. From  Fig. On the left we observe that the error in the $L^2(0,T;L^2)$-norm converges quadratically to zero for implicit and explicit* methods. For the explicit method, the convergence rate is also quadratic, however, we observe stagnation of the errors for spatial refinements $s>1$. The comparison of the results with $l=6$ and $l=7$ suggests that this stagnation tends to become less relevant when the time step size is chosen smaller.

For the errors in the  $L^2(0,T;H^1)$-norm (right plot in Fig.~\ref{Fig4}) we observe linear convergence w.r.t. spatial discretization as expected. The explicit method again underperformes and we observe a stagnation which is, however, less visible as compared to the $L^2(0,T;L^2)$ error.

\begin{figure}[t]   
	\centering
        \includegraphics[width=0.49\textwidth]{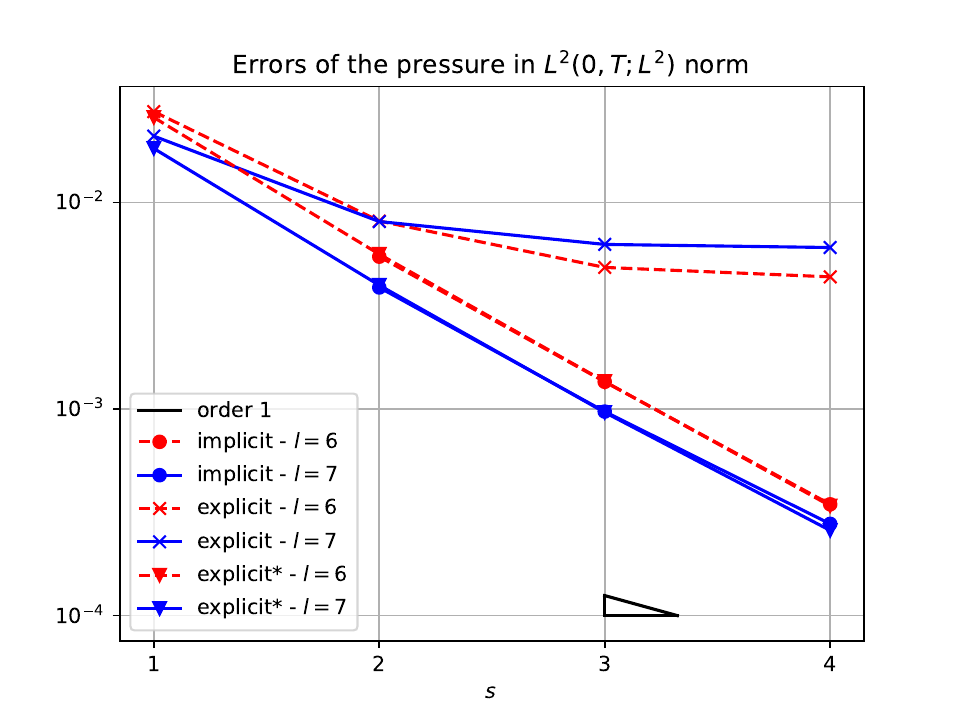}
    \caption{Errors of the pressure in the $L^2(0,T;L^2)$-norm under spatial refinement (see Table \ref{tab:L2L2_pres_h})}
    \label{Fig6}
\end{figure}

\begin{table}[h]
\centering
\begin{tabular}{c|cc|cc|cc}
\toprule
$s$ & \multicolumn{2}{c|}{implicit}   & \multicolumn{2}{c|}{explicit} & \multicolumn{2}{c}{explicit*} \\
\midrule
 & $l=0$  & $l=1$ & $l=0$  & $l=1$ &$s=0$  & $s=1$ \\
\midrule
1 & $\numprint{7.75e-4}$ & $\numprint{6.09e-4}$ & $\numprint{5.47e-4}$ & $\numprint{5.50e-4}$ & $\numprint{5.44e-4}$ & $\numprint{5.46e-4}$ \\ 
2 & $\numprint{1.37e-4}$ & $\numprint{1.37e-4}$ & $\numprint{1.34e-4}$ & $\numprint{1.43e-4}$ & $\numprint{1.30e-4}$ & $\numprint{1.31e-4}$ \\ 
3 & $\numprint{3.38e-5}$ & $\numprint{3.47e-5}$ & $\numprint{4.31e-5}$ & $\numprint{6.68e-5}$ & $\numprint{3.22e-5}$ & $\numprint{3.26e-5}$ \\ 
4 & $\numprint{9.12e-6}$ & $\numprint{1.14e-5}$ & $\numprint{2.98e-5}$ & $\numprint{6.03e-5}$ & $\numprint{8.29e-6}$ & $\numprint{9.22e-6}$ \\
\bottomrule
\end{tabular}
\caption{Errors of the predictor velocity in the $L^2(0,T;L^2)$-norm under spatial refinement}
\label{tab:L2L2_h}
\end{table}

\begin{table}[h]
\centering
\begin{tabular}{c|cc|cc|cc}
\toprule
$s$ & \multicolumn{2}{c|}{implicit}   & \multicolumn{2}{c|}{explicit} & \multicolumn{2}{c}{explicit*} \\
\midrule
 & $l=0$  & $l=1$ & $l=0$  & $l=1$ &$s=0$  & $s=1$ \\
\midrule
1 & $\numprint{2.15e-2}$ & $\numprint{2.02e-2}$ & $\numprint{1.94e-2}$ & $\numprint{1.95e-2}$ & $\numprint{1.94e-2}$ & $\numprint{1.94e-2}$ \\ 
2 & $\numprint{9.73e-3}$ & $\numprint{9.74e-3}$ & $\numprint{9.55e-3}$ & $\numprint{9.57e-3}$ & $\numprint{9.52e-3}$ & $\numprint{9.52e-3}$ \\ 
3 & $\numprint{4.80e-3}$ & $\numprint{4.80e-3}$ & $\numprint{4.78e-3}$ & $\numprint{4.89e-3}$ & $\numprint{4.73e-3}$ & $\numprint{4.73e-3}$ \\ 
4 & $\numprint{2.37e-3}$ & $\numprint{2.38e-3}$ & $\numprint{2.52e-3}$ & $\numprint{2.97e-3}$ & $\numprint{2.36e-3}$ & $\numprint{2.36e-3}$\\
\bottomrule
\end{tabular}
\caption{Errors of the predictor velocity in the $L^2(0,T;H^1)$-norm under spatial refinement}
\label{tab:L2H1_h}
\end{table}

\begin{table}[h]
\centering
\begin{tabular}{c|cc|cc|cc}
\toprule
$s$ & \multicolumn{2}{c|}{implicit}   & \multicolumn{2}{c|}{explicit} & \multicolumn{2}{c}{explicit*} \\
\midrule
 & $l=0$  & $l=1$ & $l=0$  & $l=1$ &$s=0$  & $s=1$ \\
\midrule
1 & $-$ & $-$ & $\numprint{2.76e-02}$ & $\numprint{2.10e-02}$ & $\numprint{2.58e-02}$ & $\numprint{1.83e-02}$ \\
2 & $\numprint{5.47e-03}$ & $\numprint{3.88e-03}$ & $\numprint{8.13e-03}$ & $\numprint{8.09e-03}$ & $\numprint{5.62e-03}$ & $\numprint{3.98e-03}$ \\
3 & $\numprint{1.35e-03}$ & $\numprint{9.72e-04}$ & $\numprint{4.85e-03}$ & $\numprint{6.27e-03}$ & $\numprint{1.36e-03}$ & $\numprint{9.64e-04}$ \\
4 & $\numprint{3.46e-04}$ & $\numprint{2.78e-04}$ & $\numprint{4.37e-03}$ & $\numprint{6.06e-03}$ & $\numprint{3.40e-04}$ & $\numprint{2.58e-04}$ \\
\bottomrule
\end{tabular}
\caption{Errors of the pressure field in $L^2(0,T;L^2)$-norm under spatial refinement}
\label{tab:L2L2_pres_h}
\end{table}

 \section{Conclusion}
In this paper we analysed implicit and explicit pressure correction schemes. While the implicit method has been considered in the literature before, results for the explicit variant have been missing so far. We've have extended the results and proven the expected estimates under a realistic cfl condition. 

Moreover, we presented a highly practical variant, where the explicit term is reformulated to allow representation by means of sparse matrix vector products instead of requiring numerical quadrature to set up the residual in each time step. This substantially reduces the computational effort. The theory has been extended to this case. 
We presented numerical tests which illustrates the convergence behavior of the methods.

The practical variant of the explicit method is well suited for multithreaded computations as it only consists of matrix vector multiplications except for the pressure update. In an upcoming work we plan to focus on the GPU parallelisation of the presented explicit* method. Moreover, we will extend the theoretical framework to higher order approximations.

\bibliographystyle{alpha}
\newcommand{\etalchar}[1]{$^{#1}$}

\appendix

\section{Proofs for the implicit pressure correction scheme}\label{app:implicit}

\begin{proof}[Proof to Lemma~\ref{lem:1} (Error propagation for the implicit pressure correction scheme]
The proof is close to the proof of \cite[Theorem 5.5]{Guermond1998}. We nevertheless give all details, as the explicit counterpart that is covered in the following section can be treated similarly with only few terms that require a different estimation. Furthermore, we carefully track all constants as they will enter the CFL condition of the explicit case.

Substracting \eqref{eq:pre_discrete} from \eqref{projection} gives
\begin{align}
    (\tilde\bxi_u^n - \bxi_u^{n-1}, \bchi) + k\nu(\nabla \tilde\bxi_u^{n}, \nabla \bchi)  \; =& \;k(r_h^n- p_{h}^{n-1},\div \bchi)-  k(\bbf^n,\bchi) + (\bw_h^n - \bw_h^{n-1},\bchi)  \nonumber \\ 
    &+k\nu(\nabla \bu^n, \nabla \bchi)+ k c(\tilde \bu_h^{n} ,  \tilde \bu_h^{n-1}, \bchi)\nonumber \\
     = & \; - \underbrace{k(\partial_t \bu^n - \delta_t \bu^n, \bchi)}_{=:I} - \underbrace{k( \delta_t \boleta_u^n, \bchi)}_{II} \nonumber \\
    & - \underbrace{k \Big(c(\bu^{n} ,\bu^{n}, \bchi) -  c(\tilde \bu_h^{n-1} , \tilde  \bu_h^{n}, \bchi) \Big)}_{=:III}+ \underbrace{k(j_h p^n - p_{h}^{n-1},\div \bchi)}_{=:IV} \nonumber \\
    & + \underbrace{k(r_h^n + p^n - j_h p^n,\div \bchi)}_{=:V}\label{error_equation}
\end{align}
for each $\bchi \in \bV_h$.
We choose $\bchi =\tilde\bxi_u^n $ and introduce the notation $\psi_h^n :=j_h p^n - p_h^{n-1}$. We start with the most critical term and obtain the control over pressure error using~\eqref{eq:corr} by considering
\begin{align*}
-IV & =-k(\div \tilde \bxi_u^n, \psi_h^n) =k(\div (\bxi_u^n -\tilde \bxi_u^n), \psi_h^n) = -k(\bxi_u^n -\tilde \bxi_u^n, \nabla \psi_h^n)  = -k^2(\nabla(p^n_h - p^{n-1}_h),\nabla \psi_h^n) \\
& =  -\frac{k^2}{2}\| \nabla \psi_h^n\|_0^2 + \frac{k^2}{2}\| \nabla (\psi_h^n+p_h^{n-1} - p_h^n ) \|_0^2 - \frac{k^2}{2}\| \nabla( p_h^n -p_h^{n-1}  ) \|_0^2\\
& =  -\underbrace{\frac{k^2}{2}\| \nabla \psi_h^n\|_0^2}_{=:IV_a}  + \frac{k^2}{2}\| \nabla \xi_p^n\|_0^2-  \underbrace{\frac{k^2}{2}\| \nabla( p_h^n -p_h^{n-1}  ) \|_0^2}_{=IV_b}.
\end{align*}
The second equality above is obtained by using 
\begin{align}
    (\bu-\bv,\bw) = \frac1{2} ( \| \bw\|_0^2 - \| \bw + \bv -\bu \|_0^2 + \| \bu - \bv\|_0^2 )\label{eq:err1}
\end{align}
with $\bu=\nabla p_h^n$, $\bv=\nabla p_h^{n-1}$ and $\bw=\nabla \psi_h^n$.
For  term $IV_a$ we use the fact that 
\[
\psi_h^n =j_h p^n - p_h^{n-1} (= j_h p^n - j_h p^{n-1} + j_h p^{n-1} - p_h^{n-1} )=  k \delta_t j_h p^n_h +\xi_p^{n-1}
\]
and obtain
\begin{align*} 
    IV_a= \frac{k^2}{2}\|\nabla \psi^n_h\|_0^2 \leq \frac{k^4}{2}(1+\frac1{k}) \|\nabla \delta_t j_h p^n\|_0^2 + \frac{k^2}{2}(1+k)\|\nabla \xi_p^{n-1} \|_0^2
\end{align*}
by employing $(a+b)^2\leq (1+\epsilon)a^2 + (1+\frac1{\epsilon})b^2$ for any $\epsilon>0$. The first term in $IV_a$ is estimated further as
\begin{align}
\label{grad-increment-press-projection}
\|\nabla \delta_t j_h p^n\|_0^2  &\lesssim  \|\nabla \delta_t p^n\|_0^2 \lesssim \|\partial_t p^n\|_{L^{\infty}(0,T; H^1(\Omega))}^2
\end{align}
where we used the $H^1$ stability of the projection $j_h$ and \eqref{time_derivative}.
For the term $ IV_b$, we use the stability of ~\eqref{eq:corr}, i.e.
\begin{align*}
    IV_b= \frac{k^2}{2}\| \nabla (p_h^n -p_h^{n-1}) \|_0^2 & \leq  \frac{k}{2}\| \tilde \bxi_u^n - \bxi_u^{n} \|_0 \| \nabla (p_h^n -p_h^{n-1}) \|_0 \\
    & \leq \frac{1}{4} \| \tilde \bxi_u^n - \bxi_u^{n} \|_0^2 +  \frac{k^2}{4} \|\nabla (p_h^n -p_h^{n-1}) \|_0^2 \\
    & \leq \frac{1}{2}\| \tilde \bxi_u^n - \bxi_u^{n} \|_0^2 \\
    & = \frac{1}{2}\| \tilde \bxi_u^n \|_0^2  - \frac{1}{2}\|\bxi_u^{n} \|_0^2
\end{align*}
due to $\| \bxi_u^n \|_0^2 + \| \tilde \bxi_u^n - \bxi_u^{n} \|_0^2   - \| \tilde \bxi_u^{n} \|_0^2 =0 $ which is a consequence of testing \eqref{eq:corr} with $\bxi_u^n$. For the term $V$, together with \eqref{projection-error} and \eqref{interpolation} we observe
\begin{align*}
    V & = k(r_h^n + p^n - j_h p^n,\div \bxi_u^n) = k(r_h^n,\div \bxi_u^n) - k(\nabla\boleta_p^n, \tilde \bxi_u^n) \\
    & \leq \frac{\sqrt{k}}{\sqrt{\nu}}\|r_h^n\|_0\sqrt{k\nu}\|\nabla\bxi_u^n \|_0 + k\|\boleta_p^n\|_0\|\tilde \bxi_u^n \|_0^2\\
    & \lesssim  \frac{k}{\nu}\|r_h^n\|_0^2+  c_{Y}k\nu \|\nabla\bxi_u^n \|_0^2 + k^2\|\boleta_p^n\|_0^2 + c_{Y}\|\tilde \bxi_u^n \|_0^2\\
    & \lesssim k h^{2r} (\| \bu^n\|_{r+1,2}^2 + k\| p^n\|_{r,2}^2) + c_{Y} k\nu\|\nabla\bxi_u^n \|_0^2  + c_{Y}\|\tilde \bxi_u^n \|_0^2
\end{align*}
The nonlinear term gives
\begin{align*}
    III &= k \Big(c(\bu^{n} ,\bu^{n}, \bchi) -  c( \tilde \bu_h^{n-1} , \tilde \bu_h^{n}, \bchi) \Big) \\
    &= k c(\bu^{n} - \bu^{n-1} ,\bu^{n}, \bchi) + kc(\bu^{n-1} -  \bw_h^{n-1} ,\bu^{n}, \bchi) + kc( \bw_h^{n-1}, \bu^{n}- \bw_h^{n}, \bchi)\\
    &\quad +kc(  \bw_h^{n-1} - \tilde\bu_h^{n-1},\bw_h^{n}, \bchi) +kc(  \tilde \bu_h^{n-1} ,\bw_h^{n}-\tilde \bu_h^{n}, \bchi)
\end{align*}
Each term is estimated as follows: for the first term by using H\"older's inequality we obtain
\begin{align*}
    k c(\bu^{n} - \bu^{n-1} ,\bu^{n},  \tilde \bxi_u^n) \leq \; & k^{2}\|\delta_t \bu^{n}\|_6 \|\nabla \bu^{n}\|_{3}  \| \tilde \bxi_u^n\|_0 \\
     \lesssim \;& k^{2} \|\delta_t\bu\|_{1,2} \|\bu\|_{L^{\infty}(0,T;W^{1,3}(\Omega))} \| \tilde \bxi_u^n\|_0\\
     \lesssim \;&  C_{1,n} k^4+  c_{Y}\| \tilde \bxi_u^n\|_0^2
\end{align*}
where 
\[
C_{1,n} = \|\partial_t\bu\|_{L^{\infty}(0,T;H^1(\Omega))}^2 \|\bu\|^2_{L^{\infty}(0,T;W^{1,3}(\Omega))}.
\]
From \eqref{nonlinest3} and \eqref{stab2} we have
\begin{align*}
    kc(\bu^{n-1} -  \bw_h^{n-1} ,\bu^{n},  \tilde \bxi_u^n)   \leq \;& k \| \boleta_u^{n-1} \|_{1,2} \| \bu^{n}\|_{2,2}  \| \tilde \bxi_u^n\|_0\\
    \lesssim \;& C_{2,n} k^2 h^{2r} +  c_{Y} \|\tilde \bxi_u^n\|_0^2
\end{align*}
 where $C_{2,n} = \|\bu\|^2_{L^{\infty}(0,T;H^{r+1}(\Omega))}\|\bu\|^2_{L^{\infty}(0,T;W^{2,2}(\Omega))} $. As next, from the \eqref{nonlinest2} and \eqref{stab2} we have
\begin{align*}
   kc( \bw_h^{n-1}, \bu^{n}- \bw_h^{n}, \tilde \bxi_u^n)   \leq \;& k
   \Big( \|\bw_h^{n-1}\|_{\infty} + \| \bw_h^{n-1}\|_{1,3} \Big) \| \boleta_u^{n}\|_{1,2}  \| \tilde \bxi_u^n\|_0\\
   \lesssim \;& C_{3,n} k^2 h^{2r} + c_{Y}\| \tilde \bxi_u^n\|_0^2
\end{align*}
where $C_{3,n} = \|\bu \|^2_{L^{\infty}(0,T;W^{2,2}(\Omega))}  \|\bu\|^2_{L^{\infty}(0,T;H^{r+1}(\Omega))}$. Finally, again using \eqref{stab2} combined with \eqref{nonlinest3} we have
\begin{align*}
    kc(  \tilde\bxi_h^{n-1},\bw_h^{n}, \tilde \bxi_u^n)  \leq \; & k\|\tilde \bxi_h^{n-1}\|_{1,2} \Big( \| \bw_h^n\|_{1,3} + \| \bw_h^n\|_{\infty} \Big) \| \tilde \bxi_u^{n}\|_0\\
    \leq \; & k\|\tilde \bxi_h^{n-1}\|_{1,2} \| \bu^n\|_{2,2}  \| \tilde \bxi_u^{n}\|_0\\
    \leq \; & \sqrt{\frac{C_{4,n}k}{2}} (\sqrt{ k \nu}\|\nabla  \tilde \bxi_h^{n-1}\| +  \|\tilde \bxi_h^{n-1}\|)\| \tilde \bxi_u^{n}\|_0\\
    \lesssim \; & C_{4,n}k \| \tilde \bxi_u^{n}\|_0^2 +c_Y\| \tilde \bxi_u^{n-1}\|_0^2 + c_Y \nu k \|\nabla  \tilde \bxi_h^{n-1}\|_0^2 
\end{align*}
with $C_{4,n} = 2\min \{(1 + c_p^2) \nu^{-1}, c_{inv}^2 k (1+h^{-2})\}   \|\bu^n\|_{2,2}^2$
and the last term in $III$ vanishes. This leads to 
\begin{align*}
     III \lesssim & \; C_{1,n} k^2 + ( C_{2,n}+ C_{3,n})k h^{2r} +   (C_{4,n}k + 4 c_Y) \big(\|\bxi_u^{n-1}\|_0^2 + \| \tilde \bxi_u^{n}-  \bxi_u^{n-1}\|_0^2 \big)\\
     &+ c_Y \big(\|\bxi_u^{n-2}\|_0^2 + \| \tilde \bxi_u^{n-1}-  \bxi_u^{n-2}\|_0^2\big)  + c_Y \nu k \|  \nabla \tilde \bxi_u^{n-1}\|_0^2 
\end{align*}
Moreover for the remaining terms on the right-hand-side we have from \eqref{time_approx}
\begin{align*}
    I = k(\partial_t \bu^n - \delta_t \bu^n, \tilde \bxi_u^n) \lesssim &   k^2 \|\partial_t \bu^n - \delta_t \bu^n\|_{0}^2 + c_Y\| \tilde \bxi_u^n\|_0^2  \\
    \lesssim &   k^3  \| \partial_{tt}\bu^n \|_{L^{2}(0,T;L^2(\Omega))}^2 + c_Y\| \tilde \bxi_u^n\|_0^2  
\end{align*}
and using \eqref{time_approx_boleta}
\begin{align*}
    II = k(\delta_t \boleta_u^n, \tilde \bxi_u^n) 
    \lesssim &  k^2 \|\delta_t \boleta_u^n \|_{0}^2 +  c_Y\| \tilde \bxi_u^n\|_0^2  \\
    \lesssim & k h^{2r+2} \| \partial_{t}\bu^n \|_{L^{2}(0,T;W^{r+1,2}(\Omega))}^2 +   c_Y\| \tilde \bxi_u^n\|_0^2  
\end{align*}
gives the claim.
\end{proof}

\section{Technical estimates for estimating the pressure error}
\label{app:pressure}

\begin{proof}[Proof to Lemma~\ref{approx:1:p}]
Testing the error equation \eqref{err-eq-explicit} with $\bchi = \tilde \bxi_u^1$ gives
\begin{align*}
  \frac1{2}\big(  \| \tilde \bxi_u^1 \|_0^2  +   \| \tilde \bxi_u^1 - \bxi_u^0 \|_0^2 - \| \bxi_u^{0}\|_0^2 \big) + k\nu\| \nabla \tilde \bxi_u^{1}\|^2_0\leq &  \underbrace{k\nu(\nabla (\tilde \bxi_u^{1}-\tilde \bxi_u^0),\nabla \tilde \bxi_u^{1})}_{=:I} 
    + \underbrace{k (\|\partial_t \bu^1 - \delta_t \bu^1 \|_{0} + \|\delta_t \boleta_u^1\|_{0})\| \tilde \bxi_u^{1}\|_0}_{=:II} \\
    & + \underbrace{k \Big(c(\bu^{1} ,\bu^{1}, \tilde \bxi_u^{1}) -  c(\tilde \bu_h^{0} , \tilde  \bu_h^{0}, \tilde \bxi_u^{1}) \Big)}_{=:III}\\
    & +  \underbrace{k(j_h p^0 - p_{h}^{0},\div \tilde \bxi_u^{1})}_{=0} +  \underbrace{k(r_h^1 + p^1 - j_h p^{0},\div \tilde \bxi_u^{1})}_{:=IV}
\end{align*}
Using the inverse inequality and the assumption $\| \tilde \bxi_u^{0} \|_0=0$ we obtain
\begin{align*}
    I & \leq k\nu\| \nabla (\tilde \bxi_u^{1}-\tilde \bxi_u^0)\|_0 \| \nabla \tilde \bxi_u^{1}\|_0 \leq\frac{k\nu}{h} c_{inv} \| \tilde \bxi_u^{1}-\tilde \bxi_u^0 \|  \| \nabla \tilde \bxi_u^{1}\| \\ 
    & \lesssim \frac{k\nu}{2h^2} c^2_{inv}\| \tilde \bxi_u^{1}-\tilde \bxi_u^0 \|_0^2   +  2 c_Y k \nu  \| \nabla \tilde \bxi_u^{1}\|_0^2 \leq  \frac1{2}\| \tilde \bxi_u^{1} \|_0^2  +  c_Y k \nu \| \nabla \tilde \bxi_u^{1}\|_0^2
\end{align*}
where we used $C_{\star,1}<1$ in last inequality.
The term $II$ is estimated as in the proof of Lemma \ref{lem:1}. For the term $IV$, we have
\begin{align*}
    IV&= k(r_h^1, \div \tilde \bxi_u^1) -k(\nabla (p^1 - j_h p^{0}), \tilde \bxi_u^1) \\
    & \leq  k\nu^{-1} \| r_h^1\|_0 k\nu \| \div \tilde \bxi_u^1 \|_0   + k \big(\| \nabla \boleta_p^1 \|_0 + k \| \nabla \delta_t j_h p^1 \|_0\big)\|  \tilde \bxi_u^1\|_0 \\
    & \lesssim  k^2 h^{2r} \|\bu^1\|^2_{r+1,2}   +  k^2 h^{2r} \|p\|_{r,2}^2  + k^4\| \partial_t  p \|_{L^{\infty}(0,T;H^1(\Omega))}^2 + c_Y\|  \tilde \bxi_u^1\|^2_0 +  c_Y k \nu \|  \nabla \tilde \bxi_u^1\|^2_0
  \end{align*}
where we used \eqref{grad-increment-press-projection}.
The nonlinear term gives as in \eqref{exp_nonlin_term}
\begin{align*}
    III &= k c(\bu^1 - \bu^0 ,\bu^1, \bxi_u^{1}) + kc(\bu^0 -  \bw_h^0 ,\bu^1, \bxi_u^{1}) +
    kc( \bw_h^0, \bu^1- \bw_h^1, \bxi_u^{1}) \\
    &+kc(  \bw_h^0 - \tilde\bu_h^0,\bw_h^1, \bxi_u^{1}) + kc(  \tilde \bu_h^{0} ,\bw_h^{1}-\bw_h^0,  \tilde \bxi_u^{1})+      kc(  \tilde \bu_h^{0} ,\bw_h^0 -\tilde \bu_h^{0},  \tilde \bxi_u^{1})
\end{align*}
with $c(  \bw_h^0 - \tilde\bu_h^0,\bw_h^1, \bxi_u^{1})=c(  \tilde \bu_h^{0} ,\bw_h^0 -\tilde \bu_h^{0},  \tilde \bxi_u^{1}) =0$, since $\| \tilde \bxi_u^{0} \|_0=0$. 
Remaining terms are estimated as in Lemma \ref{lem:1}:
\begin{align*}
    kc(\bu^{1} - \bu^0 ,\bu^{1}, \tilde \bxi_u^{1})  & \lesssim k^4 \| \partial_t \bu\|^2_{L^{\infty}(0,T;H^1(\Omega))} \| \bu^1\|^2_{2,2} +  c_Y\| \tilde \bxi_u^{1}\|^2_0,\\
    kc(\bu^{0} -  \bw_h^{0} ,\bu^{1}, \tilde \bxi_u^{1}) &\lesssim k^2 h^{2r}\| \bu^0\|^2_{r+1,2}  \|\bu^1\|_{2,2}^2 +  c_Y\| \tilde \bxi_u^{1}\|_0^2,\\
kc( \bw_h^{0}, \bu^{1}- \bw_h^{1}, \tilde \bxi_u^{1}) & \lesssim  k^2 h^{2r}\| \bu^0\|^2_{2,2}\| \bu^1\|^2_{r+1,2} + c_Y \| \tilde \bxi_u^{1}\|^2_0,\\
    kc(  \tilde \bu_h^{0} ,\bw_h^{1}-\tilde \bu_h^{0},  \tilde \bxi_u^{1}) & = kc(  \tilde \bu_h^{0} ,k \delta_t \bw_h^{1},  \tilde \bxi_u^{1}) =  k^2((\tilde \bu_h^{0} \cdot \nabla)\delta_t \bw_h^{1},  \tilde \bxi_u^{1}) + \frac{k^2}{2}((\div \tilde \bu_h^{0}) \delta_t \bw_h^{1},  \tilde \bxi_u^{1})\\
    &\lesssim k^2 \big(\|\tilde \bu_h^{0}  \|_{\infty} \|\nabla \delta_t \bw_h^{1}\|_0 +\|\nabla \tilde \bu_h^{0}  \|_{\infty} \| \delta_t \bw_h^{1}\|_0 \big) \| \tilde \bxi_u^{1}\|_0\\
    &= k^2 \big(\| \bw_h^{0}  \|_{\infty} \|\nabla \delta_t \bw_h^{1}\|_0 +\|\nabla  \bw_h^{0}  \|_{\infty} \| \delta_t \bw_h^{1}\|_0 \big) \| \tilde \bxi_u^{1}\|_0\\
    &\lesssim k^4  \| \bu^{0}  \|^2_{3,2}\| \partial_t \bu\|^2_{L^{\infty}(0,T;H^1(\Omega))} +  c_Y\| \tilde \bxi_u^{1}\|^2_0.
    \end{align*}
Hence, we arrive at 
\begin{align*}
\|\delta_t \tilde \bxi_u^1\|^2_0 + k \nu\|\delta_t \nabla \tilde \bxi_u^1\|^2_0  \lesssim   k^2 D_1 + h^{2r} D_2.
\end{align*}
For the pressure, we have 
\begin{align*}
\| \nabla (\xi_p^{1} - \xi_p^{0}) \|^2_0 & \lesssim \| \nabla j_h (p^{1} - p^{0}) \|^2_0  +  \| \nabla (p_h^{1} - p_h^{0}) \|^2_0 \lesssim k^2 \| \nabla \delta_t p^{1} \|^2_0  + k^{-2}\|  \tilde\bxi_u^{1} \|^2_0 \\
& \lesssim k^2 \| \partial_t p^{1} \|^2_{L^{\infty}(0,T;H^1(\Omega))}  + \|  \delta_t \tilde\bxi_u^{1} \|^2_0 + k^{-2}\|  \tilde\bxi_u^{0} \|^2_0\\
& \lesssim k^2 \| \partial_t p^{1} \|^2_{L^{\infty}(0,T;H^1(\Omega))}  +  k^2 D_1 + h^{2r} D_2
\end{align*}
\end{proof}

\begin{lemma}[Auxiliary estimate for proofing the pressure error] \label{lem:Appendix}
    For $C_{\star,n} \lesssim 1 $ sufficiently small it holds 
    \begin{align*}
        T_4^n(\delta_t\tilde\bxi_u^n) \lesssim & k^2 + h^{2r} +  c_Y\big( \| \delta_t\tilde\bxi_u^{n-1}\|_0^2 + k\nu\|\nabla \delta_t\tilde\bxi_u^{n-1}\|^2 + \| \delta_t\tilde\bxi_u^n\|_0^2 + k\nu\|\nabla \delta_t\tilde\bxi_u^n \|^2\big) \\
        & + \big(k \nu^{-1} \| \tilde\bxi_u^{n-1}\|_0^2 + k^2  \|\nabla  \tilde\bxi_u^{n-1}\|_0^2\big) \|\partial_t \bu\|_{L^{\infty}(0,T;H^2(\Omega))}^2
        \end{align*}
\end{lemma}
\begin{proof}
Note that
\begin{align*}
c(\bu^{n} ,\bu^{n}, \delta_t\tilde\bxi_u^n) - c(\bu^{n-1} ,\bu^{n-1}, \delta_t\tilde\bxi_u^n)  & = c(\bu^{n} - \bu^{n-1} ,\bu^{n-1}, \delta_t\tilde\bxi_u^n) + c(\bu^{n} ,\bu^{n} - \bu^{n-1}, \delta_t\tilde\bxi_u^n)\\
& = c(\bu^{n} - \bu^{n-1} ,\bu^{n-1} - \bw_h^{n-1}, \delta_t\tilde\bxi_u^n)+  c(\bu^{n} - \bu^{n-1} ,\bw_h^{n-1}, \delta_t\tilde\bxi_u^n) \\
& + c(\bu^{n} - \bw_h^n ,\bu^{n} - \bu^{n-1}, \delta_t\tilde\bxi_u^n) + c(\bw_h^n ,\bu^{n} - \bu^{n-1}, \delta_t\tilde\bxi_u^n)\\
&=: I + II + III + IV.
\end{align*}
We further split each term 
\begin{align*}
II =& c(\bu^{n} - \bu^{n-1} ,\bw_h^{n-1}, \delta_t\tilde\bxi_u^n) = \underbrace{c(\bu^{n} - \bu^{n-1} - (\bw_h^n-\bw_h^{n-1}) ,\bw_h^{n-1}, \delta_t\tilde\bxi_u^n)}_{=:II_a} + c(\bw_h^n-\bw_h^{n-1} ,\bw_h^{n-1}, \delta_t\tilde\bxi_u^n)\\
=& II_a + \underbrace{c(\bw_h^n-2\bw_h^{n-1} + \bw_h^{n-2} ,\bw_h^{n-1}, \delta_t\tilde\bxi_u^n)}_{=:II_b}+ c(\bw_h^{n-1} - \bw_h^{n-2} ,\bw_h^{n-1}, \delta_t\tilde\bxi_u^n)\\
=& II_a + II_b+ \underbrace{c(\bw_h^{n-1} - \bw_h^{n-2} ,\bw_h^{n-1}-\tilde \bu_h^{n-1}, \delta_t\tilde\bxi_u^n)}_{II_c} + c(\bw_h^{n-1} - \bw_h^{n-2} ,\tilde \bu_h^{n-1}, \delta_t\tilde\bxi_u^n)\\
=& II_a + II_b+ II_c + \underbrace{c(\bw_h^{n-1} - \bw_h^{n-2}-(\tilde \bu_h^{n-1} - \tilde \bu_h^{n-2}) ,\tilde \bu_h^{n-1}, \delta_t\tilde\bxi_u^n)}_{=:II_d} + c(\tilde \bu_h^{n-1} - \tilde \bu_h^{n-2} ,\tilde \bu_h^{n-1}, \delta_t\tilde\bxi_u^n)\\
=& II_a + II_b+ II_c + II_d + c(\tilde \bu_h^{n-2} ,\tilde \bu_h^{n-2} - \tilde \bu_h^{n-1} , \delta_t\tilde\bxi_u^n) - c(\tilde \bu_h^{n-2} , \tilde \bu_h^{n-2}, \delta_t\tilde\bxi_u^n)+c(\tilde \bu_h^{n-1} , \tilde \bu_h^{n-1}, \delta_t\tilde\bxi_u^n)\\
=& II_a + II_b+ II_c + II_d +  \underbrace{c(\tilde \bu_h^{n-2} , \tilde \bu_h^{n-2} - \tilde \bu_h^{n-1} - (\bw_h^{n-2} - \bw_h^{n-1})  , \delta_t\tilde\bxi_u^n) }_{:=II_e} + \underbrace{c(\tilde \bu_h^{n-2} ,  \bw_h^{n-2} - \bw_h^{n-1}  , \delta_t\tilde\bxi_u^n)}_{=:II_f} \\
&+ c(\tilde \bu_h^{n-2} , \tilde \bu_h^{n-2}, \delta_t\tilde\bxi_u^n)+c(\tilde \bu_h^{n-1} , \tilde \bu_h^{n-1}, \delta_t\tilde\bxi_u^n)\\
\end{align*}  
and
\begin{align*}
    IV =& c(\bw_h^n ,\bu^{n} - \bu^{n-1}, \delta_t\tilde\bxi_u^n) =  \underbrace{c(\bw_h^n ,\bu^{n} - \bu^{n-1} - (\bw_h^{n} - \bw_h^{n-1}), \delta_t\tilde\bxi_u^n)}_{=:IV_a}  + c(\bw_h^n ,\bw_h^{n} - \bw_h^{n-1}, \delta_t\tilde\bxi_u^n) \\
    = &IV_a  + \underbrace{c(\bw_h^n - \bw_h^{n-1} ,\bw_h^{n} - \bw_h^{n-1}, \delta_t\tilde\bxi_u^n)}_{=:IV_b}   +c(\bw_h^{n-1} ,\bw_h^{n} - \bw_h^{n-1}, \delta_t\tilde\bxi_u^n) \\
    = &IV_a  + IV_b   +\underbrace{c(\bw_h^{n-1} - \tilde \bu_h^{n-1} ,\bw_h^{n} - \bw_h^{n-1}, \delta_t\tilde\bxi_u^n)}_{=:IV_c} + c(\tilde \bu_h^{n-1} ,\bw_h^{n} - \bw_h^{n-1}, \delta_t\tilde\bxi_u^n)\\
    = &IV_a  + IV_b   + IV_c + \underbrace{c(\tilde \bu_h^{n-1} ,\bw_h^{n} - 2\bw_h^{n-1} + \bw_h^{n-2}, \delta_t\tilde\bxi_u^n)}_{=:IV_d} + \underbrace{c(\tilde \bu_h^{n-1} ,\bw_h^{n-1} - \bw_h^{n-2}, \delta_t\tilde\bxi_u^n)}_{=:IV_e} 
    \end{align*}  
Therefore, we have    
\begin{align*}
-T_4^n(\delta_t\tilde\bxi_u^n) &= c(\bu^{n} ,\bu^{n}, \delta_t\tilde\bxi_u^n) - c(\bu^{n-1} ,\bu^{n-1}, \delta_t\tilde\bxi_u^n) -  c(\tilde \bu_h^{n-1} , \tilde  \bu_h^{n-1}, \delta_t\tilde\bxi_u^n) + c(\tilde \bu_h^{n-2} , \tilde  \bu_h^{n-2}, \delta_t\tilde\bxi_u^n) \\
& = I + II_a + II_b+ II_c + II_d + II_e+ II_f+ III+ IV_a  + IV_b  + IV_c  + IV_d + IV_e
 \end{align*}
Classical arguments yield to following term-by-term estimates
 \begin{align*}
 I & = c(k\delta_t\bu^{n} ,\boleta_u^{n-1}, \delta_t\tilde\bxi_u^n) \lesssim k^2\|\delta_t \bu^{n}\|^2_{\infty} \|\nabla \boleta_u^{n-1}\|^2_{0} + c_Y  \|  \delta_t \tilde \bxi_u^n\|^2_0 \\
& \lesssim  k^{2}h^{2r}\|\delta_t \bu^{n}\|^2_{2,2} \| \bu^{n-1}\|^2_{r+1,2}+ c_Y  \| \delta_t \tilde \bxi_u^n\|^2_0 
\end{align*}
\begin{align*}
II_a & = c(k\delta_t\boleta_u^{n} ,\bw_h^{n-1}, \delta_t\tilde\bxi_u^n) \leq k \|\nabla \delta_t \boleta_u^{n}\|_0 \big(\| \bw_h^n\|_{\infty} + \|  \bw_h^n\|_{1,3} \big)\| \delta_t\tilde\bxi_u^n\|_0\\
 &\lesssim k^2 \|\nabla \delta_t \boleta_u^{n}\|^2_0 \| \bu^n\|_{2,2}^2 + c_Y \| \delta_t\tilde\bxi_u^n\|^2_0
\end{align*}
\begin{align*}
II_b & = c(k^2\delta_{tt}\bw_h^{n} ,\bw_h^{n-1}, \delta_t\tilde\bxi_u^n) \leq k^2 \|\delta_{tt}\bw_h^{n}\|_{1,2} \big(\| \bw_h^{n-1}\|_{\infty} + \| \bw_h^{n-1}\|_{1,3} \big)\| \delta_t\tilde\bxi_u^n\|_0\\
&\lesssim k^4 \|\delta_{tt}\bu^{n}\|_{1,2}^2  \| \bu^{n-1}\|^2_{2,2} + c_Y\| \delta_t\tilde\bxi_u^n\|_0^2
\end{align*}
\begin{align*}
II_c & = c(k \delta_t \bw_h^{n-1},\tilde\bxi_u^{n-1}, \delta_t\tilde\bxi_u^n)    \leq k \| \tilde\bxi_u^{n-1}\|_0 \big(\|\delta_t \bw_h^{n-1}\|_{\infty} + \|\delta_t \bw_h^{n-1}\|_{1,3}\big)\| \nabla \delta_t\tilde\bxi_u^n\|_0\\
& \lesssim k \nu^{-1} \| \tilde\bxi_u^{n-1}\|_0^2 \|\delta_t \bu^{n-1}\|_{2,2}^2 + c_Y k \nu\| \nabla \delta_t\tilde\bxi_u^n\|_0^2\\
\end{align*}
Choosing $k \big(\frac{\nu}{h^2}c^2_{inv} + \frac{\|\tilde \bu_h^{n-1}\|^2_{\infty}}{\nu}\big) \leq \sqrt{c_Y}$ gives
\begin{align}
\label{estimate_complicated}
\begin{split}
II_d & = c(k\delta_t \tilde \bxi_u^{n-1} ,\tilde \bu_h^{n-1}, \delta_t\tilde\bxi_u^n) 
\lesssim  k\| \delta_t \tilde \bxi_u^{n-1}\|_0 \| \tilde \bu_h^{n-1}\|_{\infty}\| \nabla \delta_t\tilde\bxi_u^n\|_0 + \| \div \delta_t \tilde \bxi_u^{n-1}\|_0 \| \tilde \bu_h^{n-1}\|_{\infty} \| \delta_t\tilde\bxi_u^n\|_0 \\
& \lesssim  k h^{-1} c_{inv}\| \tilde \bu_h^{n-1}\|_{\infty}\| \delta_t \tilde \bxi_u^{n-1}\|_0 \| \delta_t\tilde\bxi_u^n\|_0\\
& \lesssim k^2 h^{-2}c^2_{inv} \| \tilde \bu_h^{n-1}\|^2_{\infty}\big(\| \delta_t \tilde \bxi_u^{n-1}\|^2_0 + \| \delta_t \tilde \bxi_u^{n}\|_0^2\big)\\
& \lesssim k^2 h^{-2}c^2_{inv} \nu \sqrt{c_Y} \big(\| \delta_t \tilde \bxi_u^{n-1}\|_0^2 + \| \delta_t \tilde \bxi_u^{n}\|_0^2\big)\\
& \lesssim c_Y \big(\| \delta_t \tilde \bxi_u^{n-1}\|_0^2 + \| \delta_t \tilde \bxi_u^{n}\|_0^2 \big)
\end{split}
\end{align}
Similar to the estimate \eqref{estimate_complicated}
\begin{align*}
 |II_e| & = |c(\tilde \bu_h^{n-2} , k\delta_t\tilde \bxi_u^{n-1}, \delta_t\tilde\bxi_u^n) |\\
& \lesssim k \|\tilde \bu_h^{n-2} \|_{\infty} \big(\| \delta_t \tilde \bxi_u^{n-1}\|_0 \| \nabla \delta_t\tilde\bxi_u^n\|+ \| \nabla \delta_t \tilde \bxi_u^{n-1}\|_0+ \| \delta_t\tilde\bxi_u^n\|^0\big)\\
& \lesssim k \|\tilde \bu_h^{n-2} \|^2_{\infty} \nu^{-1} \big(\| \delta_t \tilde \bxi_u^{n-1}\|_0^2 + \| \delta_t\tilde\bxi_u^n\|^2_0\big) + c_Y k\nu \| \nabla \delta_t\tilde\bxi_u^n\|^2 + c_Y k\nu \| \nabla \delta_t \tilde \bxi_u^{n-1}\|^2_0\\
& \lesssim \sqrt{c_Y} \big(\| \delta_t \tilde \bxi_u^{n-1}\|_0^2 + \| \delta_t\tilde\bxi_u^n\|^2_0\big) + c_Y k\nu \| \nabla \delta_t\tilde\bxi_u^n\|^2 + c_Y k\nu \| \nabla \delta_t \tilde \bxi_u^{n-1}\|^2_0
\end{align*}
\begin{align*}
 III & = c(\boleta_u^{n} ,k\delta_t\bu^{n}, \delta_t\tilde\bxi_u^n) 
  \lesssim k \|\nabla \boleta_u^{n}\|_0\big(\| \delta_t\bu^{n}\|_{\infty} + \|  \delta_t\bu^{n}\|_{1,3} \big)\| \delta_t\tilde\bxi_u^n\|_0\\
 & \lesssim k^2 \|\nabla \boleta_u^{n}\|^2_0 \| \delta_t\bu^{n}\|_{2,2}  + c_Y\| \delta_t\tilde\bxi_u^n\|^2_0\\
 & \lesssim k^2 h^{2r}\| \bu^n \|_{2,2}  \| \delta_t\bu^{n}\|_{2,2}  + c_Y \| \delta_t\tilde\bxi_u^n\|^2_0
\end{align*}
\begin{align*}
IV_a & = c(\bw_h^n ,k\delta_t\boleta_u^{n}, \delta_t\tilde\bxi_u^n)
\lesssim k \big(\| \nabla \delta_t\boleta_u^{n}\|_0 \|\bw_h^n\|_{\infty} + \|\div \bw_h^n\|_{3}\| \delta_t\boleta_u^{n}\|_6 \big)\| \delta_t\tilde\bxi_u^n\|_0\\
& \lesssim k^2  \| \delta_t\boleta_u^{n}\|^2_{1,2} \|\bu^n\|_{2,2}^2 + \|  \delta_t\tilde\bxi_u^n\|_0^2\\
& \lesssim k  h^{2r}\| \partial_t\bu\|^2_{L^2(0,T;W^{r+1,2}(\Omega))} \|\bu^n\|_{2,2}^2 + \|  \delta_t\tilde\bxi_u^n\|_0^2
\end{align*}
\begin{align*}
IV_b & = c(k\delta_t\bw_h^n ,k\delta_t\bw_h^n, \delta_t\tilde\bxi_u^n)
\lesssim k^2 \| \delta_t\bw_h^n\|_{1,2} \big(\|\delta_t\bw_h^n\|_{\infty} + \| \delta_t\bw_h^n\|_{1,3} \big)\|  \delta_t\tilde\bxi_u^n\|_0\\
& \lesssim k^4 \| \delta_t\bu^n\|^2_{1,2} \|\delta_t\bu^n\|_{2,2} + c_Y\|  \delta_t\tilde\bxi_u^n\|_0^2
\end{align*}
\begin{align*}
IV_c & = c(\tilde \bxi_u^{n-1} ,k\delta_t\bw_h^{n}, \delta_t\tilde\bxi_u^n) \lesssim k \| \tilde \bxi_h^{n-1}\|_{0} (\| \delta_t\bw_h^{n}\|_{\infty}  + \|\delta_t\bw_h^{n}\|_{1,3}\big)\|\nabla \delta_t\tilde\bxi_u^n\|_0\\
& \lesssim k \nu^{-1} \| \tilde \bxi_h^{n-1}\|^2_{0}\| \delta_t\bu^{n}\|_{2,2}^2 + k \nu c_Y \|\nabla \delta_t\tilde\bxi_u^n\|_0^2
\end{align*}
\begin{align*}
IV_d & = c(\tilde \bu_h^{n-1} ,k^2\delta_t\bw_h^{n}, \delta_t\tilde\bxi_u^n)
 \lesssim k^2 \| \tilde \bu_h^{n-1}\|_{\infty} \| \delta_{tt}\bw_h^{n} \|_0 \| \nabla \delta_t\tilde\bxi_u^n\|_0 + k^2 \| \tilde \bu_h^{n-1}\|_{\infty} \| \nabla \delta_{tt}\bw_h^{n} \|_0 \|  \delta_t\tilde\bxi_u^n\|_0\\
& \lesssim k^3  \| \tilde \bu_h^{n-1}\|^2_{\infty} \nu^{-1} \big( \| \delta_{tt}\bw_h^{n} \|^2_0 +\| \delta_{t}\tilde\bxi_u^n \|_0^2 \big) + c_Y k\nu \| \nabla \delta_t\tilde\bxi_u^n\|_0^2 + c_Y k\nu  \| \nabla \delta_{tt}\bw_h^{n} \|_0
\end{align*}            

\begin{align*}
    |II_f + IV_e| & = |c(\tilde \bu_h^{n-1} - \tilde \bu_h^{n-2} ,k\delta_t\bw_h^{n-1}, \delta_t\tilde\bxi_u^n)|  = |c(k\delta_t \bw_h^{n-1} ,k\delta_t\bw_h^{n-1}, \delta_t\tilde\bxi_u^n) -   c(k\delta_t \tilde \bxi_u^{n}  ,k\delta_t\bw_h^{n-1}, \delta_t\tilde\bxi_u^n)|\\
    & \lesssim k^4 \| \delta_t \bu^{n-1}\|^2_{1,2} \| \delta_{t}\bu^{n-1} \|^2_{2,2}+ c_Y \| \delta_t\tilde\bxi_u^n\|^2_0 + k^2 \|\nabla \tilde\bxi_u^n\|_0 \|\delta_t\bu^{n-1}\|_{2,2} + k \nu c_Y \|\nabla \delta_t\tilde\bxi_u^n\|_0^2
    \end{align*}            
\end{proof}    
\end{document}